\documentclass[final]{siamltex}
\usepackage[utf8]{inputenc}
\usepackage{amsmath,amssymb,amsfonts,etoolbox}
\usepackage[hidelinks]{hyperref}
\hypersetup{colorlinks=true,citecolor=blue,linkcolor=blue,
            filecolor=blue,urlcolor=blue}
\usepackage{graphicx}
\usepackage{cite}
\usepackage[font={small}]{caption}
\usepackage{url}
\usepackage{booktabs}
\usepackage{framed}
\usepackage{mathtools}
\usepackage{pgf,tikz}
\usepackage{mathrsfs}
\usepackage{mathtools, enumerate}
\usepackage{multirow}
\usetikzlibrary{arrows}

\def\dfrac{\displaystyle\frac}
\def\dsum{\displaystyle\sum}

\newcommand{\snorm}[1]{\left| #1 \right|}
\newcommand{\norm}[1]{\left\| #1 \right\|}
\newcommand{\lrbrack}[1]{\left( #1 \right)}

\newcommand{\lrBbrack}[1]{\left\{ #1 \right\}}

\newtheorem{example}{Example}[section]
\newtheorem{remark}{Remark}[section]
\newtheorem{hypothesis}{Hypothesis}[section]

\begin{document}

\title{
Fast linear barycentric rational interpolation for singular functions via scaled transformations
\thanks{This work was supported by the National Natural Science Foundation of China (No. 11771454).
The first author is partly supported by the Fundamental Research Funds for the Central Universities of Central South University (No. 2020zzts031).}
}

\author{
Desong Kong \thanks{School of Mathematics and Statistics, Central South University,
 Changsha, Hunan 410083, People's Republic of China
  ({\tt desongkong@csu.edu.cn}, {\tt xiangsh@csu.edu.cn}).}
\and Shuhuang xiang \footnotemark[2]
\thanks{Corresponding author.}
}

\maketitle

\begin{abstract}
  In this paper, applied  strictly monotonic increasing scaled maps, a kind of well-conditioned linear barycentric rational  interpolations are proposed to approximate  functions of singularities at the origin, such as $x^\alpha$ for $\alpha \in (0,1)$ and $\log(x)$. It just takes $O(N)$ flops and can achieve fast convergence rates with the choice the scaled parameter,  where $N$ is the maximum degree of the denominator and numerator. The construction of the rational interpolant couples rational polynomials in the barycentric form of second kind with the transformed Jacobi-Gauss-Lobatto points. Numerical experiments are considered which illustrate the accuracy and efficiency of the algorithms. The convergence of the rational interpolation is also considered.
\end{abstract}



\pagestyle{myheadings}
\thispagestyle{plain}
\markboth{
\uppercase{Desong Kong and Shuhuang Xiang}}
{\uppercase{
Rational Interpolation for Singular Functions
}}

\section{Introduction}\label{sec:int}
Function approximation is a classical topic in scientific computation. It is known that approximation by polynomials for analytic functions can be characterized perfectly by Bernstein's theorem \cite[Chapter 8]{Trebook2019} with exponential convergence rates. While for functions with endpoint or interior algebraic singularity, just algebraically convergent rates can be achieved (see \cite[Chapter 7]{Trebook2019} and \cite{Xiang2016,Xiang2020}).

Compared with polynomial approximation, rational approximation may achieve fast convergence rates even for functions of limited regularity. A fundamental result for rational approximation owns to Newman's work for $f(x)=|x|$  \cite{Newman1964}, who showed that the rational function
\begin{align}\label{newman}
  r(x)=x\frac{p(x)-p(-x)}{p(x)+p(-x)},\quad p(x)=\prod_{k=0}^{N-1} (x+\xi^k),\quad \xi=\exp(-\sqrt{N}),\quad x\in [-1,1]
\end{align}
can achieve \textit{root-exponential convergence}
\begin{align*}
  \frac12 e^{-9\sqrt{N}} \leq \|f-r\|_{\infty} \leq 3e^{-\sqrt{N}},\quad N \geq 4,
\end{align*}
 which is essentially much better to the one order convergence with polynomial approximation $\|f-p_N^*\|_{\infty}=O(N^{-1})$ \cite{Bernstein1914,Poussin1908}, where $p_N^*$ is the best approximation polynomial of degree $N$.

Newman's investigation triggers a whole series of contributions to improve the error estimate. Stahl's series of theoretical investigation finally answered the best uniform rational approximation for $|x|$ in \cite{Stahl1993} and also gave an extended result for $x^\alpha$ in \cite{Stahl2003}:
\begin{align}\label{eq:newmann1}
  \lim_{N \to \infty} e^{2\pi \sqrt{\alpha N}} E_{NN}(x^\alpha, [0,1]) = 4^{1+\alpha} |\sin \pi \alpha|
\end{align}
or equivalently
\begin{align}\label{eq:newmann2}
  \lim_{N \to \infty} e^{\pi \sqrt{\alpha N}} E_{NN}(|x|^\alpha, [-1,1]) = 4^{1+\alpha/2} |\sin \frac12 \pi \alpha|
\end{align}
for each $\alpha>0$. In the special case $\alpha=1$, the best uniform rational approximation bound for $|x|$ can be derived.

It is worth noting that, for rational interpolation with respect to polynomial interpolation, its theoretical analysis and algorithm implementation become more complicated.
For instance, the traditional problem of rational interpolation
is known to present two main difficulties \cite{Berrut1988}: one is the occurrence of unattainable points, and the other is the appearance of poles in the interval of interpolation.

 To avoid the poles in the interpolation interval, extremely well conditioned rational interpolants, based upon the barycentric formula of second kind \cite{Dupuy1948} to approximate the target function $f$
 \begin{align}\label{eq:rat}
  r(x)=\left. \dsum\limits_{i=0}^N \dfrac{\lambda_i f(x_i) }{x-x_i} \right/ \dsum\limits_{i=0}^N \dfrac{\lambda_i}{x-x_i},
\end{align}
are  proposed in  Berrut et.al. \cite{Balten1999,Berrut1988,Berrut2014} by choosing simple weights  $\lambda_i = (-1)^i$, or simplified Chebyshev weights \cite{Salzer1972}
\begin{align}\label{eq:weig}
  \lambda_i =(-1)^i \delta_i,\quad \text{for} \quad
  \delta_i =
  \begin{dcases}
    1, & i=1,2\ldots,N-1,\\
    1/2, & i=0, N.
  \end{dcases}
\end{align}
The interpolation with weights \eqref{eq:weig} converges exponentially when the $x_i$ are the (shifted) Chebyshev points of the second kind if $f$ is analytic in a Bernstein ellipse enclosing the interpolation interval, as it  coincides with the
polynomial interpolant \cite{Berrut2014}. This property is conserved with any conformal map of such nodes \cite{Balten1999}.
Such conformal maps have also been extensively studied recently by Kosloff and
 Tal-Ezer \cite{Kos1993}, Bayliss and Turkel \cite{Bayliss1992}, Tee and Trefethen \cite{Tee2006} and Hale and Tee \cite{Hale2009} etc.

Based on Chebfun  system,   Deun and Trefethen \cite{Deun2011}  presented a robust implementation of the Carath\'{e}odory-Fej\'{e}r (CF) method  for rational approximation. CF approximation can be seen as a perfectly variable alternative to the best rational approximation for smooth functions \cite{Trefethen1981,Trefethen1983}. However, CF approximations are far easier to  be implemented.

More recently, two adaptive algorithms named \texttt{aaa} and \texttt{minimax}
are proposed by Nakatsukasa, Sete and Trefethen \cite{Nakat2018}, and    Filip, Nakatsukasa,  Trefethen  and Beckermann \cite{Filip2018}, respectively. These two algorithms are both built on the barycentric representation of rational functions in a third fashion
\begin{align*}
  r(x)=\left. \dsum\limits_{i=0}^N \dfrac{\alpha_i }{x-x_i} \right/ \dsum\limits_{i=0}^N \dfrac{\beta_i}{x-x_i}.
\end{align*}
The  \texttt{aaa} algorithm offers a speed,  flexible and robust  implementation with complexity $O(MN^3)$ flops, where $M$ is the number of the sample set.

The \texttt{minimax} is developed by making use of rational barycentric representations whose support
points are chosen in an adaptive fashion.
A similarly adaptive approach is also established to study the rational minimax approximation of complex functions on arbitrary domains \cite{Nakat2020a}. To obtain the best rational approximation however is still problematic and generally is an NP-hard problem. In many applications, it is not yet necessary to approximate a function with the best one.

Abundant works about rational approximation
for  solving  PDEs with singular solutions have also been developed. Trefethen et al.  \cite{Gopal2019,Trefethen2020} introduced a root-exponential approximation by rational functions with the poles preassigned clustering exponentially near the singularity, called ``linghtning'' method.

 Motivated by the series of works about rational approximation for functions of singularity,
in this paper we are interested in the rational interpolation for functions of singularity at the origin.

Suppose $f(x)$ is defined on $\Omega:=[0,T]$ and has a singularity at $x=0$ ($f(x) \sim x^\alpha$
with $\alpha \in (0,1)$ or $f(x) \sim \log(x)$). The linear barycentric rational interpolation is represented by \eqref{eq:rat}
with barycentric weights $\{ \lambda_i \}_{i=0}^N$ associated with distinct interpolation nodes $\{ x_i \}_{i=0}^N$.

It is of particular importance to note that from the root-exponential convergence rate \eqref{eq:newmann1} about  the best uniform rational approximation for   $x^\alpha$ in \cite{Stahl1993,Stahl2003}, to find a ``good" conformal map $g$ to get fast convergence is impossible even though the compound function$f(g)$ is analytic. If so, the rational approximation \eqref{eq:rat} is exponentially convergent contradicted with  Stahl's result \eqref{eq:newmann1}.

Here we  introduce  a strictly monotonic increasing scaled map $g: [-1,1] \to \Omega$:
\begin{align}\label{eq:refpts}
  x = g(y) := T \lrbrack{ \frac{y+1}{2} }^{s/\alpha},\quad y \in [-1,1]
\end{align}
for some integer $s$ and $\alpha \in (0,1)$ reflecting the singularity of $f$ near the origin.
With map \eqref{eq:refpts}, \eqref{eq:rat} is also a linear rational interpolation of type $(N,N)$ satisfying $r(x_i) = f(x_i)=f(g(y_i))$ for $x_i=g(y_i)$ and $y_i \in [-1,1]$ provided $\lambda_i \neq 0$, $i=0,1,\ldots,N$. Particularly, in the case $s=\alpha$, \eqref{eq:rat} degenerates to a polynomial interpolation. 
In general, we can replace $s/\alpha$ by $s$ in \eqref{eq:refpts} for $s \in (1, \infty)$:
\begin{align}\label{eq:refpts2}
  x = g(y) := T \lrbrack{ \frac{y+1}{2} }^{s},\quad y \in [-1,1].
\end{align}

\begin{figure}[pt]
  \centering
  \includegraphics[width=.6\textwidth]{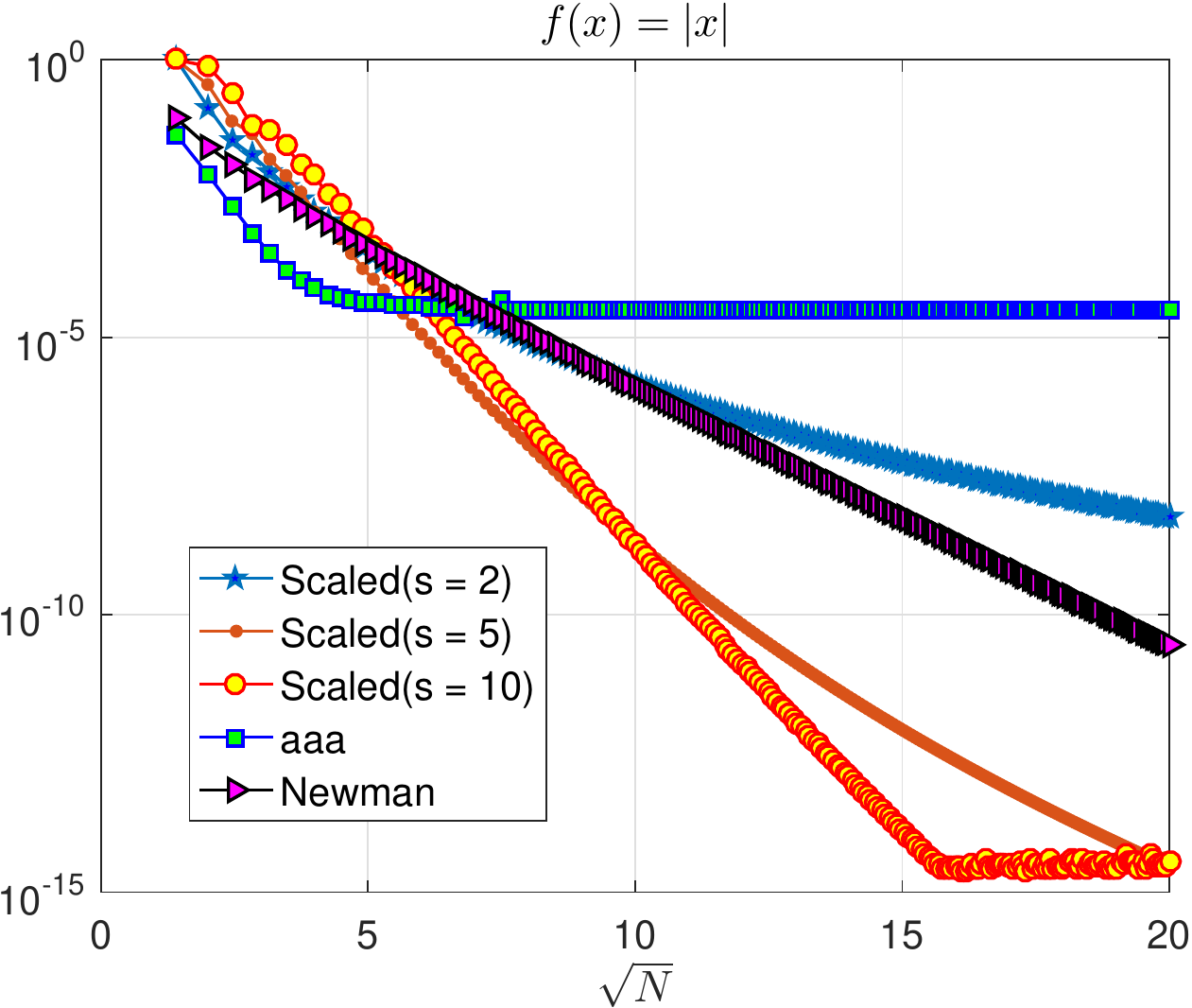}
  \caption{Comparison of three kinds of methods for approximation of $f(x) = |x|$ on $[-1,1]$: Newman's method \eqref{newman}, AAA method (denoted by `aaa' and implemented by \texttt{aaa(@abs, chebpts(2e4), 'degree', N)}) and  the proposed method (denoted by `Scaled' and implemented by calling MATLAB code \texttt{ratscale} proposed in Section~\ref{sec:num}) with various values $s=2, 5, 10$ for even values of $N$ from 2 to 400 in maximum error in \texttt{xx} (defined in Section~\ref{sec:num}).}
  \label{fig:comp}
\end{figure}

From \eqref{eq:rat}, it takes only $O(N)$ operations  for evaluating the interpolating function $r$  with given weights $\lambda_i$ and interpolation nodes $x_i$ by \eqref{eq:refpts} or \eqref{eq:refpts2}.  If we choose weights \eqref{eq:weig} coupled with transformation \eqref{eq:refpts2},  compared with Newman's method  \eqref{newman} and AAA method, the proposed rational interpolant \eqref{eq:rat}  for approximation of $f(x)=|x|$,   converges faster than the other two  for large $N$, while  Newman's method  \eqref{newman} takes $O(N^2)$ flops and AAA $O(MN^3)$. In addition, the root-exponential may  be achieved for $s = 10$ (see {\sc Fig.}~\ref{fig:comp}).

\begin{figure}[pt]
  \centering
  \begin{minipage}[c]{0.315\hsize}
    \centering
    \includegraphics[width=1\textwidth]{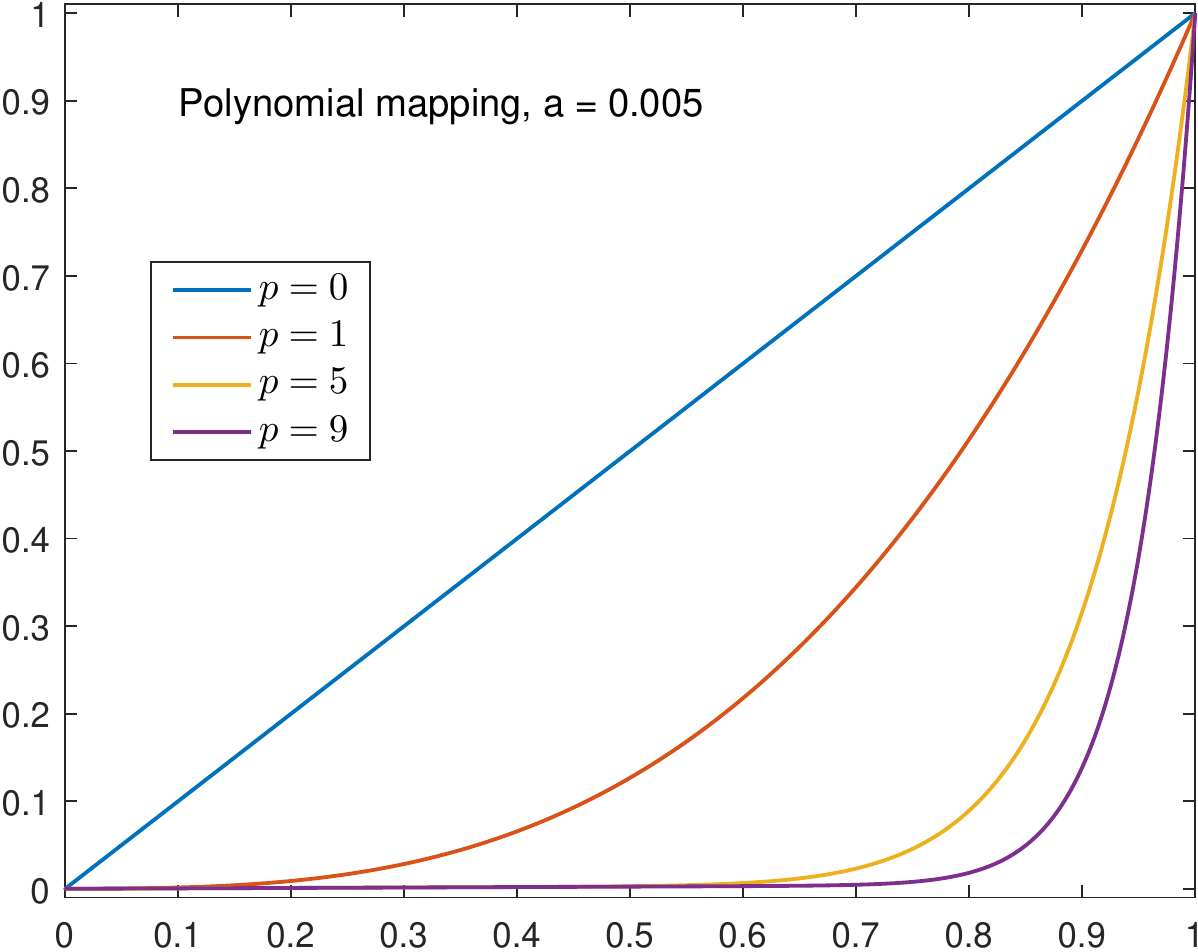}
  \end{minipage}
  \begin{minipage}[c]{0.315\hsize}
    \centering
    \includegraphics[width=1\textwidth]{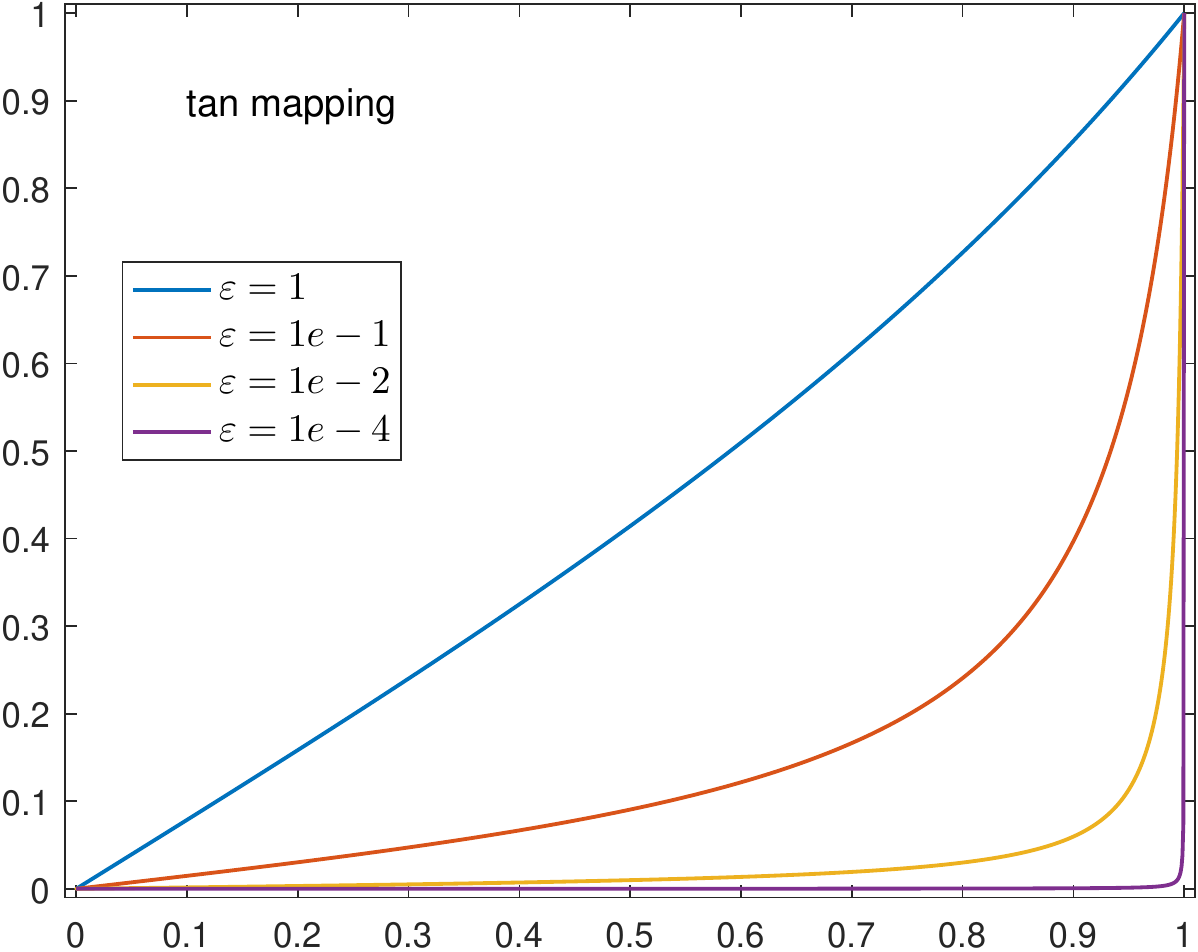}
  \end{minipage}
  \begin{minipage}[c]{0.315\hsize}
    \centering
    \includegraphics[width=1\textwidth]{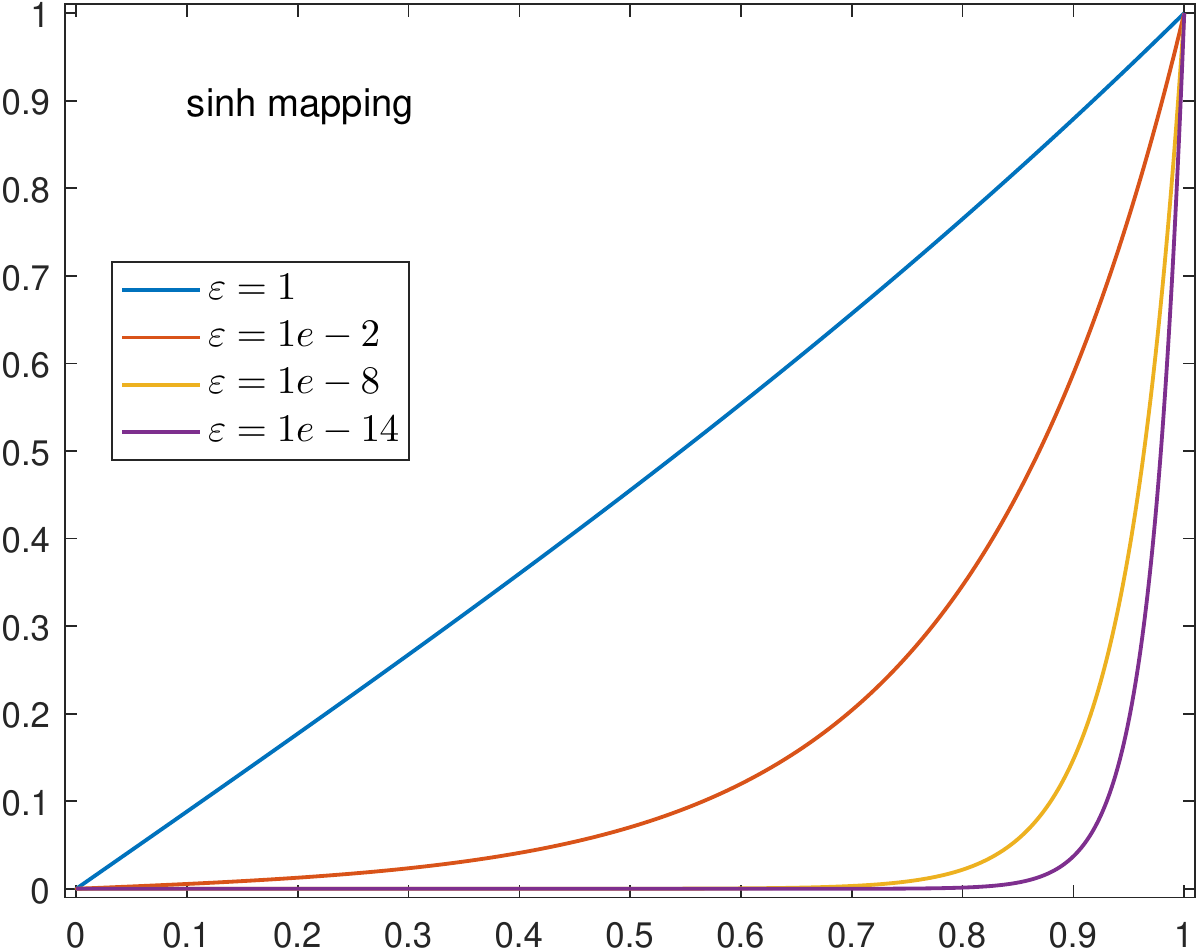}
  \end{minipage} \\
  \begin{minipage}[c]{0.317\hsize}
    \centering
    \includegraphics[width=1\textwidth]{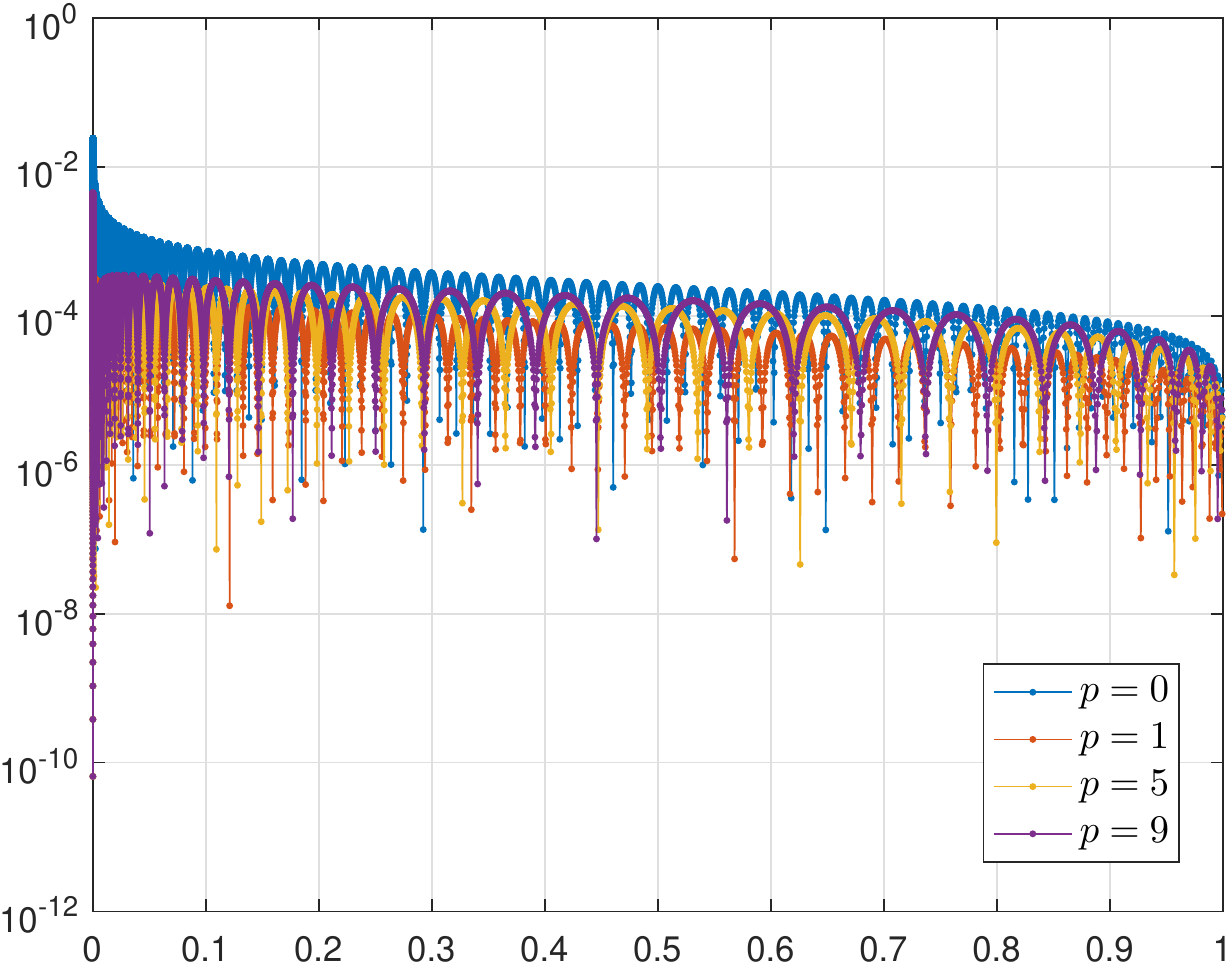}
  \end{minipage}
  \begin{minipage}[c]{0.317\hsize}
    \centering
    \includegraphics[width=1\textwidth]{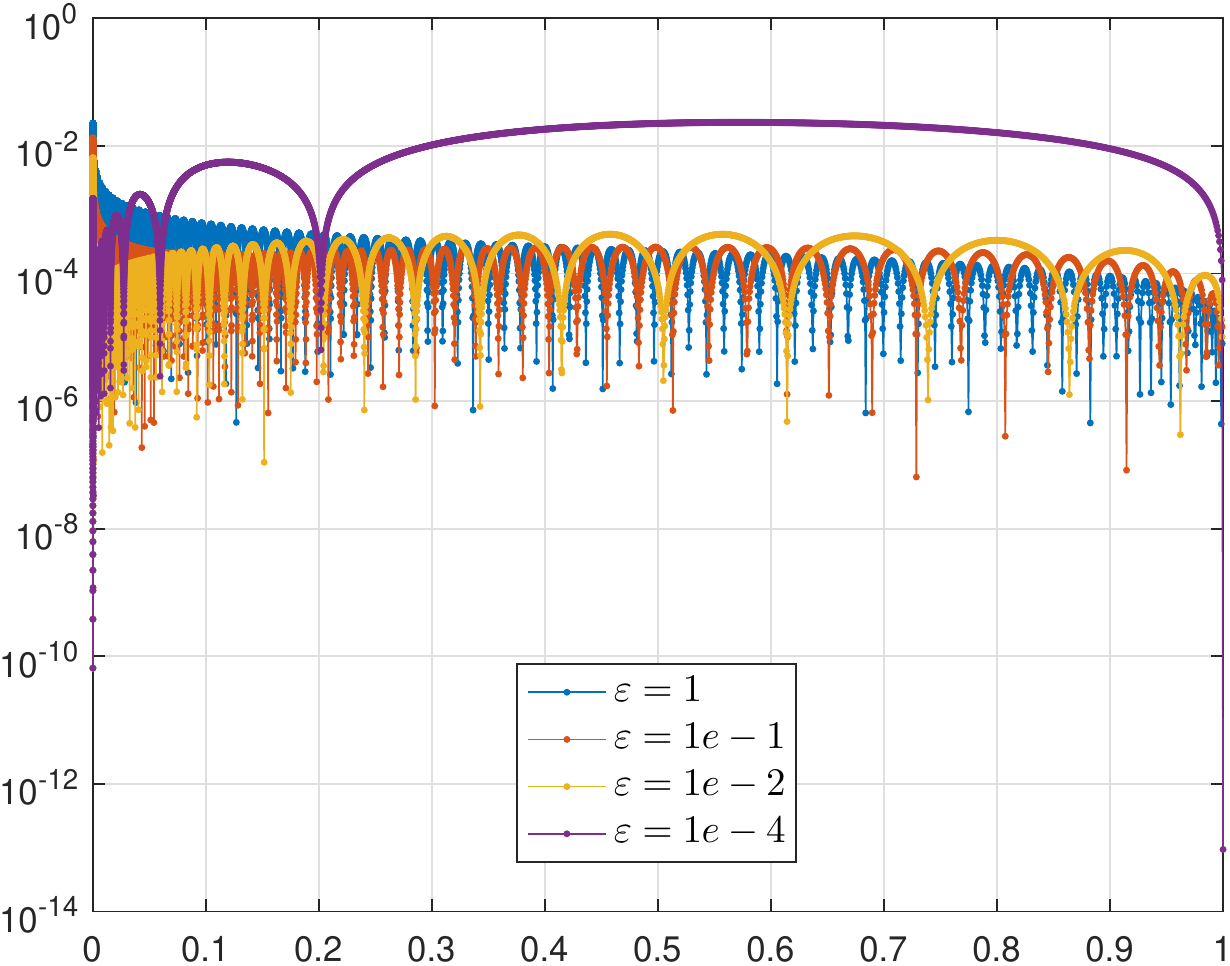}
  \end{minipage}
  \begin{minipage}[c]{0.317\hsize}
    \centering
    \includegraphics[width=1\textwidth]{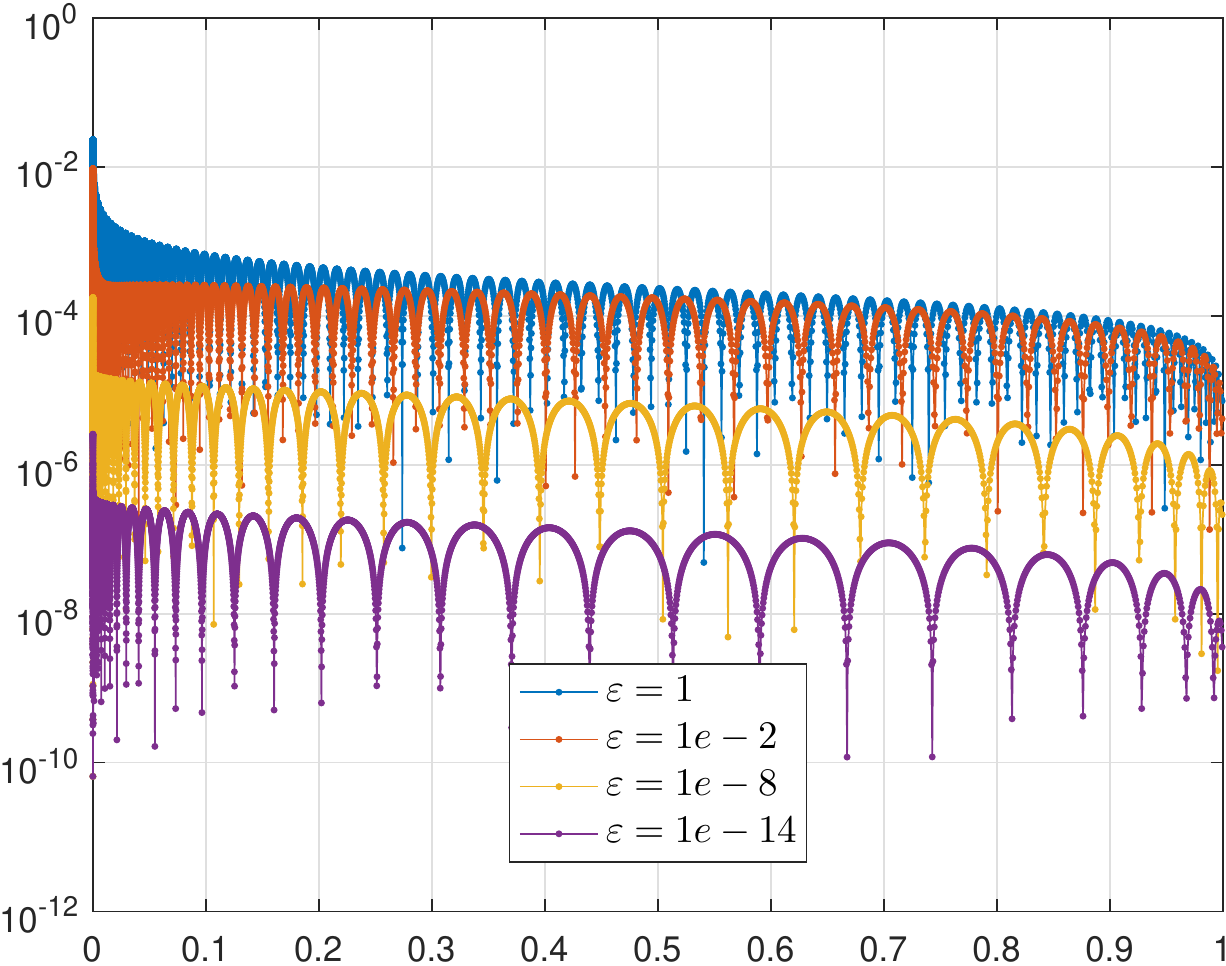}
  \end{minipage} \\
  \begin{minipage}[c]{0.315\hsize}
    \centering
    \includegraphics[width=1\textwidth]{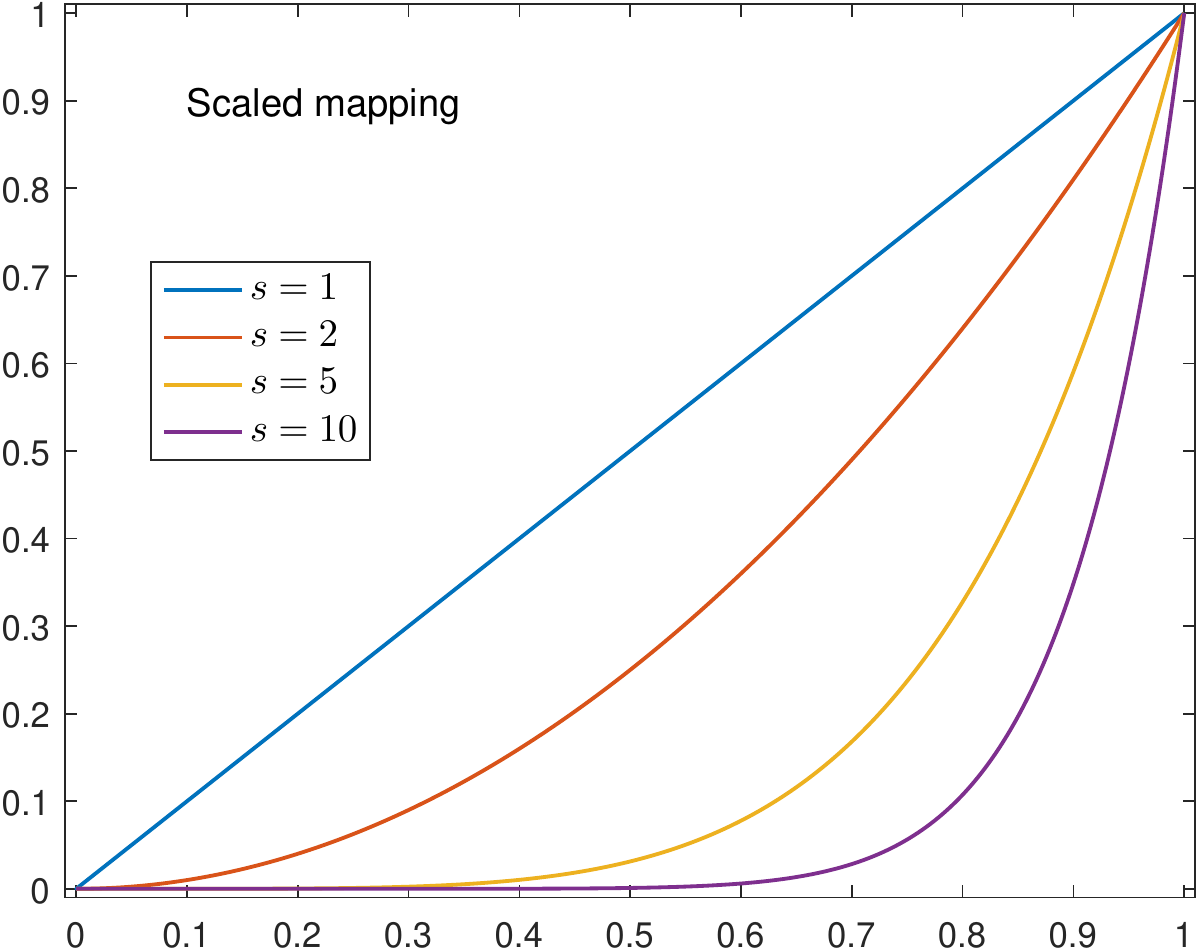}
  \end{minipage}
  \begin{minipage}[c]{0.315\hsize}
    \centering
    \includegraphics[width=1\textwidth]{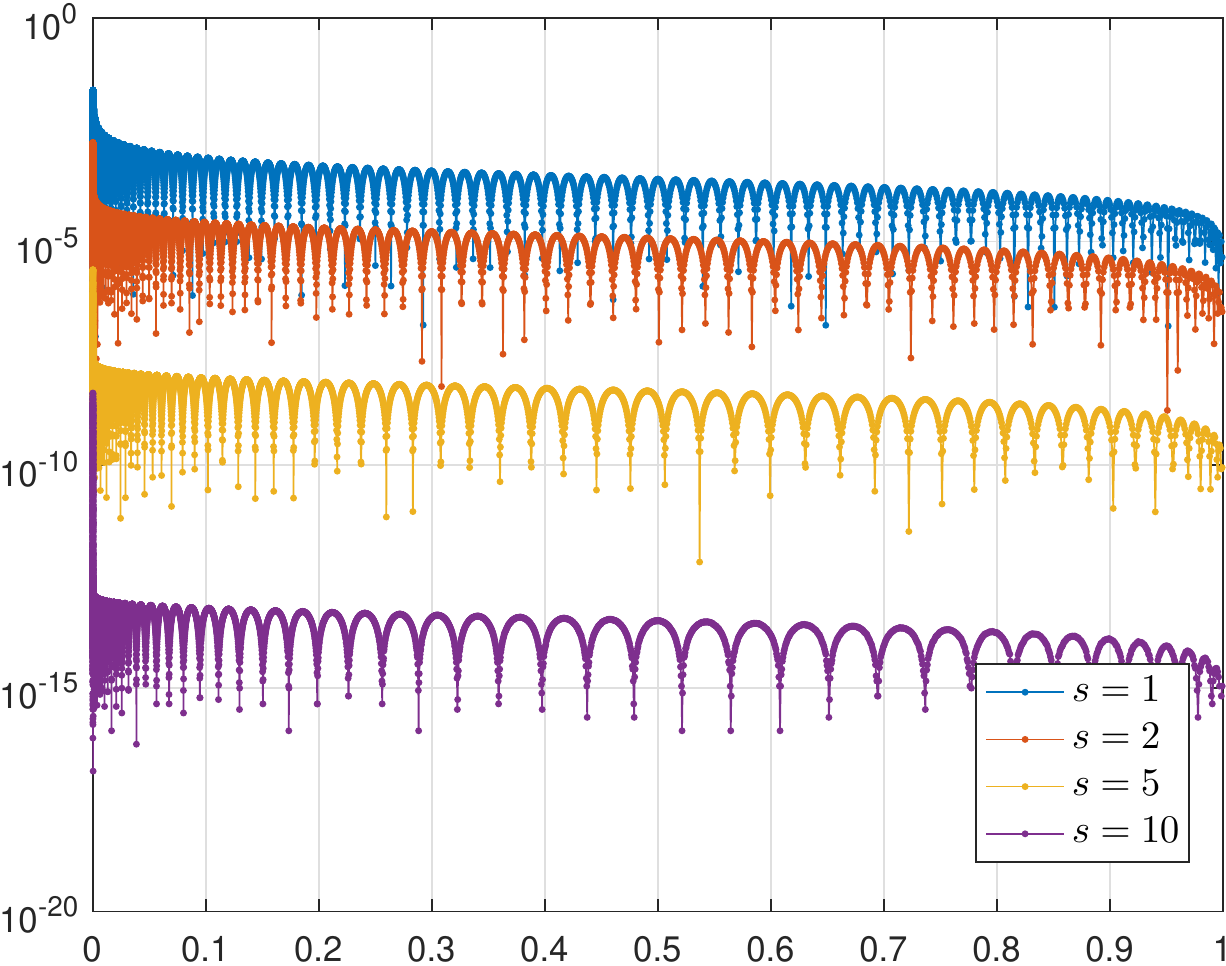}
  \end{minipage}
  \begin{minipage}[c]{0.315\hsize}
    \centering
    \includegraphics[width=1\textwidth]{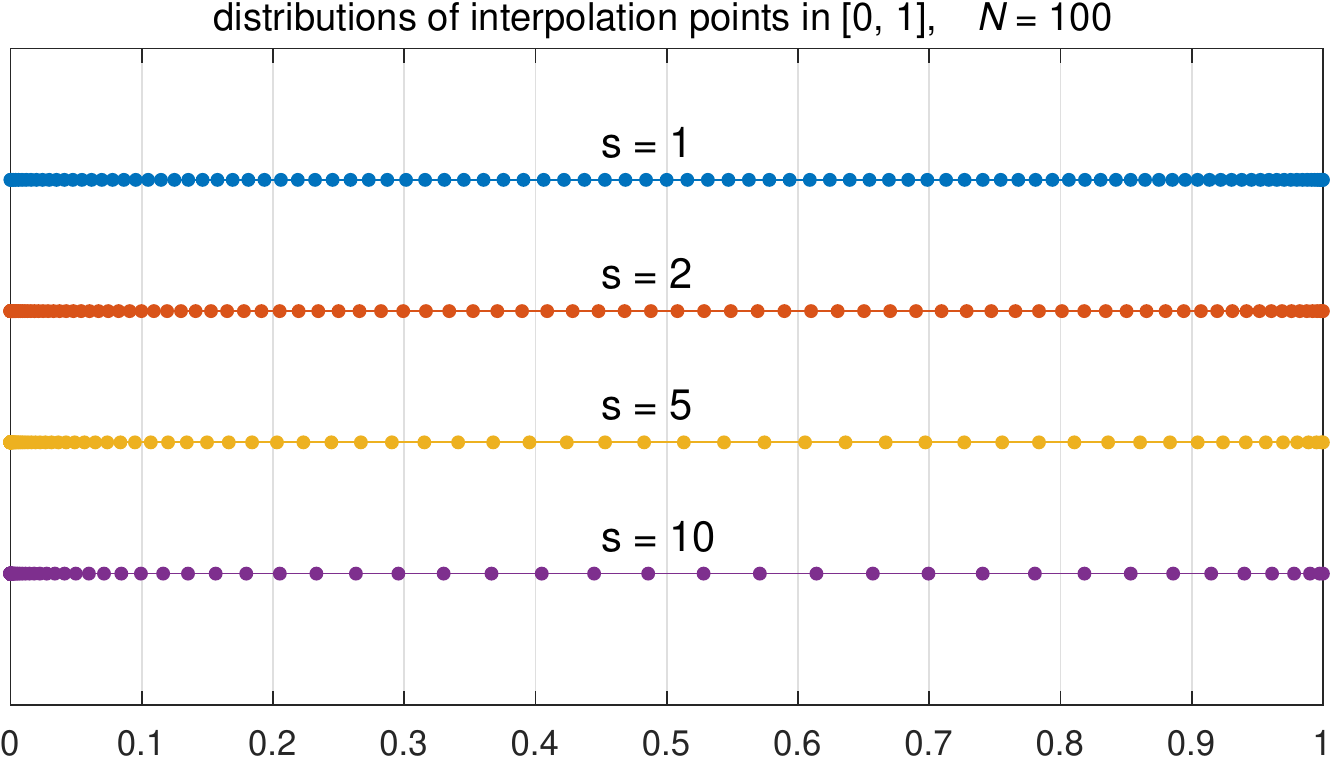}
  \end{minipage}
  \caption{Absolute error $ |f(x)-r(x)|$ corresponding to the  rational interpolation \eqref{eq:rat} with the conformal maps (second  row) given by \eqref{comps1}-\eqref{comps3} (first row), and the proposed rational approximation with \eqref{eq:refpts2} (third row) for $f(x) = x^{1/\pi}$ with $N=100$. In each plots we simulate the rational interpolations with various parameters. In the last column of the third row, the distributions of
nodes are displayed.}
  \label{fig:comp2}
\end{figure}

 We further compare with the rational interpolation \cite{Balten1999} with weights \eqref{eq:weig} coupled with the following three conformal maps for $f(x) = x^\alpha$: polynomial mapping \cite{Alexand2009}, tan-mapping \cite{Bayliss1992} and sinh-mapping \cite{Tee2006}
\begin{subequations}
\begin{align}
  x &= (ay+y^{2p+1})/(1+a),\quad a>0,\quad \mbox{and} \quad p \in \mathbb{N}; \label{comps1} \\
  x &= \varepsilon \tan(y\tan^{-1}(1/\varepsilon)),\quad \varepsilon>0; \label{comps2} \\
  x &= \varepsilon \sinh(y\sinh^{-1}(1/\varepsilon)),\quad \varepsilon>0. \label{comps3}
\end{align}
\end{subequations}
All these maps included the scaled \eqref{eq:refpts} transform the nodes concentrated at the singular point $x=0$. Numerical results  in {\sc Fig.}~\ref{fig:comp2}  show that the proposed rational approximation  performs much better than  the rational interpolation \cite{Balten1999} with weights \eqref{eq:weig} and  conformal maps \eqref{comps1}-\eqref{comps3}.

For functions of logarithmic singularity, we  introduce the following monotonic increasing map
\begin{align}\label{eq:refpts_log}
  x = g(y) = e^y \quad \mbox{for} \quad  y \in [\log(x_0), \log(T)]
\end{align}
for $x \in [x_0, T] \subset \Omega$ with $x_0>0$. The implementation of rational interpolation $r$ associated with \eqref{eq:refpts_log} is quite simple and shown in {\sc Fig.}~\ref{fig:ratlog}. Numerical results in Section~\ref{sec:exten} illustrate that the rational interpolant \eqref{eq:rat} with map \eqref{eq:refpts_log}  can achieve exponential convergence.

\begin{figure}[hbtp]
\begin{verbatim}
  function rat = ratlog(f, dom, N)
  % Input:    f: function handle
  %           N: # of interpolation nodes
  %         dom: domain of f
  % Output: rat: rational approximation to f

  gy = @(x) exp(x);                           % map
  xmin = log(min(dom)); xmax = log(max(dom)); % reference domain
  [x, ~, uk] = chebpts(N ,[xmin,xmax]);       % wts {uk}
  xi = gy(x);                                 % pts {xi}
  rat = @(x) bary(x, f(xi), xi, uk);          % rational interpolation
\end{verbatim}
\caption{MATLAB code \texttt{ratlog} for  the linear barycentric rational interpolation with scaled points for logarithmic singular functions.}
\label{fig:ratlog}
\end{figure}

The rest of this paper is outlined as follows: In Section~\ref{sec:pre}, we recall some preliminaries about the properties of barycentric rational interpolation \eqref{eq:rat}. In Section~\ref{sec:num}, various numerical experiments are simulated with the proposed rational functions to approximate  functions  of algebraic singularity. 
An extensional barycentric rational interpolant to functions of logarithmic singularity is presented in Section~\ref{sec:exten}.  Convergence  is considered in Section~\ref{sec:conv}. A brief discussion is included in the final section.

\section{Preliminaries}\label{sec:pre}
Recall that a rational function is of type $(m, n)$  \cite{Stahl1993}  if it can be written in the form of $p/q$, where $p$ and $q$ are both polynomials of degrees less than or equal to $m$ and $n$, respectively.

Let us firstly cite a theorem which proves that formula \eqref{eq:rat} is indeed an interpolation function.
\begin{theorem}[rational barycentric representation {\cite[Theorem 2.1]{Nakat2018}}]
  Let $x_0$, $\ldots$, $x_m$ be an arbitrary set of distinct complex numbers. As $f_0, \ldots, f_m$ range over all complex values and $\lambda_0, \ldots, \lambda_m$ range over all nonzero complex values, the functions $r(x)$ given by \eqref{eq:rat}
  range over the set of all rational functions of type $(m, m)$ that have no poles at the points $x_i$. Moreover, $r(x_i)=f_i$ for each $i$.
\end{theorem}

In \eqref{eq:rat}, there are two sets needed to be determined, i.e., $\{ \lambda_i \}_{i=0}^N$ and $\{ x_i \}_{i=0}^N$. For $\lambda_i = (-1)^i$, it induces  an extremely well conditioned rational interpolation with no poles \cite{Berrut1988}.
 The only major drawback, however, is the slow convergence that it is just order one with respect to the maximum length of step size.

To obtain a higher order of the approximation,
Floater and Hormann \cite{Floater2007} considered a rational function by blending local approximations to form a global one.  For equispaced points,
the Lebesgue constant of this rational interpolation grows logarithmically with $N$, but exponentially with $d$ (degree of the blended polynomials), if $f$ is analytic  \cite{Bos2011}.
The choice of the optimal $d$ for a given finite $N$ is cleared up in \cite{Guettel2012}.

Another approaches to obtain a well-conditioned rational interpolation is utilizing the simplified Chebyshev weights \eqref{eq:weig}.
The idea of barycentric rational formula associated with these transformed points is widely applied. In \cite{Balten1999}, a conformal map $g$ is introduced and generalized in \cite{Tee2006} to enlarge the ellipse of analyticity of $f$.
To reduce ill-conditioning of the Chebyshev differential matrices in solving differential equations with large gradients, some kinds of conformal map are proposed recently (see \cite{Tee2006,Jafari2014} and references therein).

The barycentric rational interpolation  can also be applied to approximate the solution of differential equations \cite{Balten2001, Berrut2001}. In this purpose, we need the following proposition about differential matrix.
\begin{proposition}[\cite{Balten1999,Sch1986}]\label{pro:der}
  Let $r$ be a rational function given in \eqref{eq:rat} with $\lambda_i \neq 0, i=0, \ldots, N$. Assume that $\xi$ is not a pole of $r$; then for $\xi \neq x_i, i=0,1,\ldots,N$,
  \begin{align*}
    \left. \frac{r^{(m)}(\xi)}{m !}=\sum_{i=0}^{N} \frac{\lambda_i}{\xi-x_i} r[(\xi)^{m}, x_i] \right/ \sum_{i=0}^{N} \frac{\lambda_i}{\xi-x_i}, \quad m \geq 0,
  \end{align*}
  and
  \begin{align*}
    \frac{r^{(m)}(x_{i})}{m !} = \left. - \lrbrack{ \sum_{i=0, i \neq j}^{N} \lambda_i r[(x_j)^{m}, x_i] } \right/ \lambda_{j}, \quad 0 \leq j \leq N, \quad m \geq 1.
  \end{align*}
\end{proposition}
The notation $(\xi)^{m}$ is used here to indicate the $m$-fold argument $\xi, \xi, \ldots, \xi,$ and $r[x_{0}, \ldots, x_{k}]$ denotes the $k$-th order divided difference of $r$ with
\begin{align*}
  r[\underbrace{\xi, \ldots, \xi}_{m+1}]=\frac{r^{(m)}(\xi)}{m !}.
\end{align*}
Making use of the above formula, we can compute the first derivative of the function $r(x)$ at the interpolation points $x_{i}$ by constructing the differentiation matrices $D^{(1)}$ whose entries are given by
\begin{align}\label{eq:diffm1}
  D_{i j}^{(1)}=
  \begin{dcases}
    \frac{\lambda_j}{\lambda_i} \frac{1}{x_i-x_j}, & i \neq j, \\
    -\sum_{k=0, k \neq i}^{N} \frac{\lambda_k}{\lambda_i} \frac{1}{x_i-x_k}, & i=j.
  \end{dcases}
\end{align}
The implementation of the first differential matrix in MATLAB is constructed by calling \texttt{bcamatrix} in \cite{Tee2006}.

\section{Fast linear barycentric rational interpolation}\label{sec:num}
Assume $f(x) \sim x^\alpha$  with $\alpha \in (0,1)$ defined  on $\Omega=[0,T]$.  In \eqref{eq:rat}, appealing to \eqref{eq:weig} for $\lambda_i$ and choosing
\begin{align}\label{eq:refpts1}
  x_i= T \lrbrack{ \frac{\cos(i\pi/N)+1}{2} }^{s/\alpha},\quad i=0,\ldots,N,
\end{align}
yield the rational approximation to $f(x)$.
The third row of {\sc Fig.}~\ref{fig:comp2} shows the scaled transformation $g(y)$ and distributions of interpolation points \eqref{eq:refpts1} in $[0,1]$ with different values of $s$ ($\alpha=1$). One can observe that as  $s$ enlarges, more points are accumulating at $x=0$, which is reasonable for approximating functions of  a singular point  at the origin.

The MATLAB code \texttt{ratscale} is shown in {\sc Fig.}~\ref{fig:ratscale} to implement the rational interpolation $r$ associated with \eqref{eq:weig} and \eqref{eq:refpts1}. The function command \texttt{bary} we used represents the barycentric rational formula \eqref{eq:rat}, and its implementation is included in Chebfun and can also be found in \cite{Berrut2004}.
With the available code, one can try any functions to test its efficiency and robustness.

\begin{figure}[hptb]
\begin{verbatim}
  function rat = ratscale(f, N, dom, s, alp)
  % Input:    f: function handle
  %           N: # of interpolation nodes
  %         dom: domain of f
  %           s: a positive integer
  %         alp: singularity of function f
  %    An optimal 5th argument specifies singularity of f near zero.
  %    If omitted, alp = 1, and s can be token as any positive number.
  % Output: rat: rational approximation to f

  if nargin<5, alp = 1; end;                  % default value alp
  gy = @(x) max(dom)*x.^(s/alp);              % map
  [x, ~, uk] = chebpts(N, [0,1]); xi = gy(x); % pts {xi} & wts {uk}
  rat = @(x) bary(x, f(xi), xi, uk);          % rational interpolation
\end{verbatim}
\caption{MATLAB code \texttt{ratscale} for the linear barycentric rational interpolation with scaled points.}
\label{fig:ratscale}
\end{figure}

The MATLAB program \texttt{ratscale} is just making use of approximation of a function defined on $\Omega$, specifically on $[0,1]$.  If one want to approximate a function defined on $[-T,T]$, for example,
$f(x) \sim |x|^\alpha$ for $x \in [-T, T]$ with  a singularity at the origin, setting
  \begin{align*}
    g(y) = \begin{dcases}
      T \lrbrack{\frac{1+y}{2}}^{s/\alpha}, & x \in [0,T], \\
      -T \lrbrack{\frac{1+y}{2}}^{s/\alpha}, & x \in [-T, 0),
    \end{dcases}
  \end{align*}
an alternative approach is replacing
the third line in {\sc Fig.}~\ref{fig:ratscale} by
\begin{verbatim}
  x = chebpts(N+1, [0,1]); xi = [-gy(x(end:-1:2)); gy(x(2:end))];
  uk = (-1).^(1:length(xi)); uk(1) = uk(1)/2; uk(end) = uk(end)/2;
\end{verbatim}

\vspace{0.36cm}
We consider various applications with \texttt{ratscale} to illustrate the efficiency and accuracy in the rest of this section. To measure the accuracy, the discrete infinite norm are used at points
\begin{verbatim}
  xx = linspace(0,1,10000); xx = xx.^8;
\end{verbatim}
for functions defined on $[0,1]$, or
\begin{verbatim}
  xx = linspace(0,1,10000); xx = xx.^8;  xx = [-xx(end:-1:2), xx];
\end{verbatim}
for functions defined on $[-1,1]$. Both cluster at the origin.

\begin{example}[Absolute value function]\label{ex:abs}\rm
  As the prototype of approximation of nonsmooth functions, absolute value function, $|x|$ on $[-1,1]$, is one of two famous problems related to rational approximation \cite[Chapter 25]{Trebook2019}. Executing the following codes in MATLAB, less than 0.002sec is needed to obtain the interpolant with error 5.58e-5 in double-precision arithmetic.
\begin{verbatim}
  dom = [-1,1]; s = 2; N = 40; f = @(x) abs(x);
  tic, rat = ratscale(f, N/2, dom, s); toc
  err = norm(f(xx)-rat(xx), inf)
  plot(xx, rat(xx)-f(xx))
\end{verbatim}
\begin{figure}[pt]
  \centering
  \includegraphics[width=.45\textwidth]{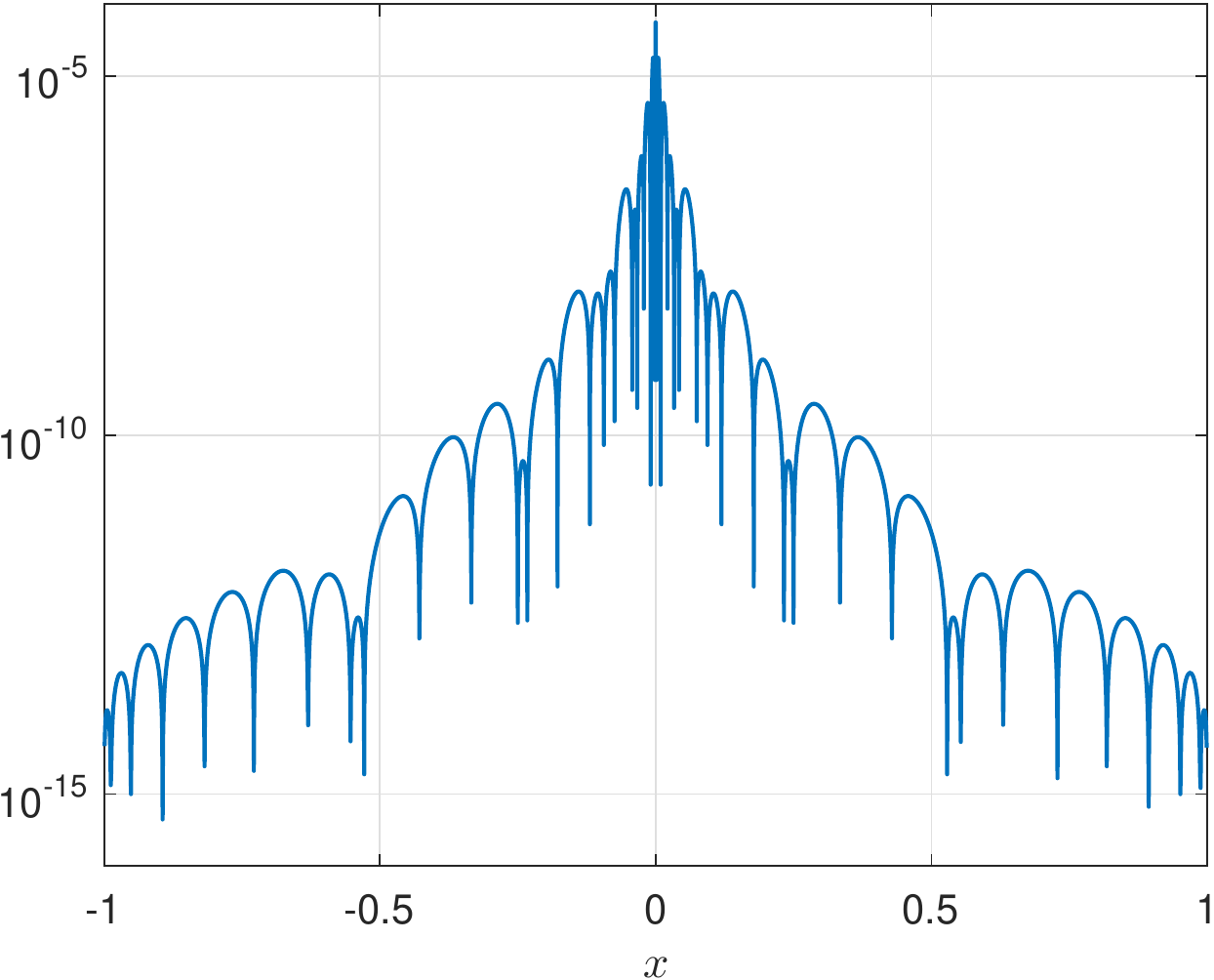}
  \includegraphics[width=.45\textwidth]{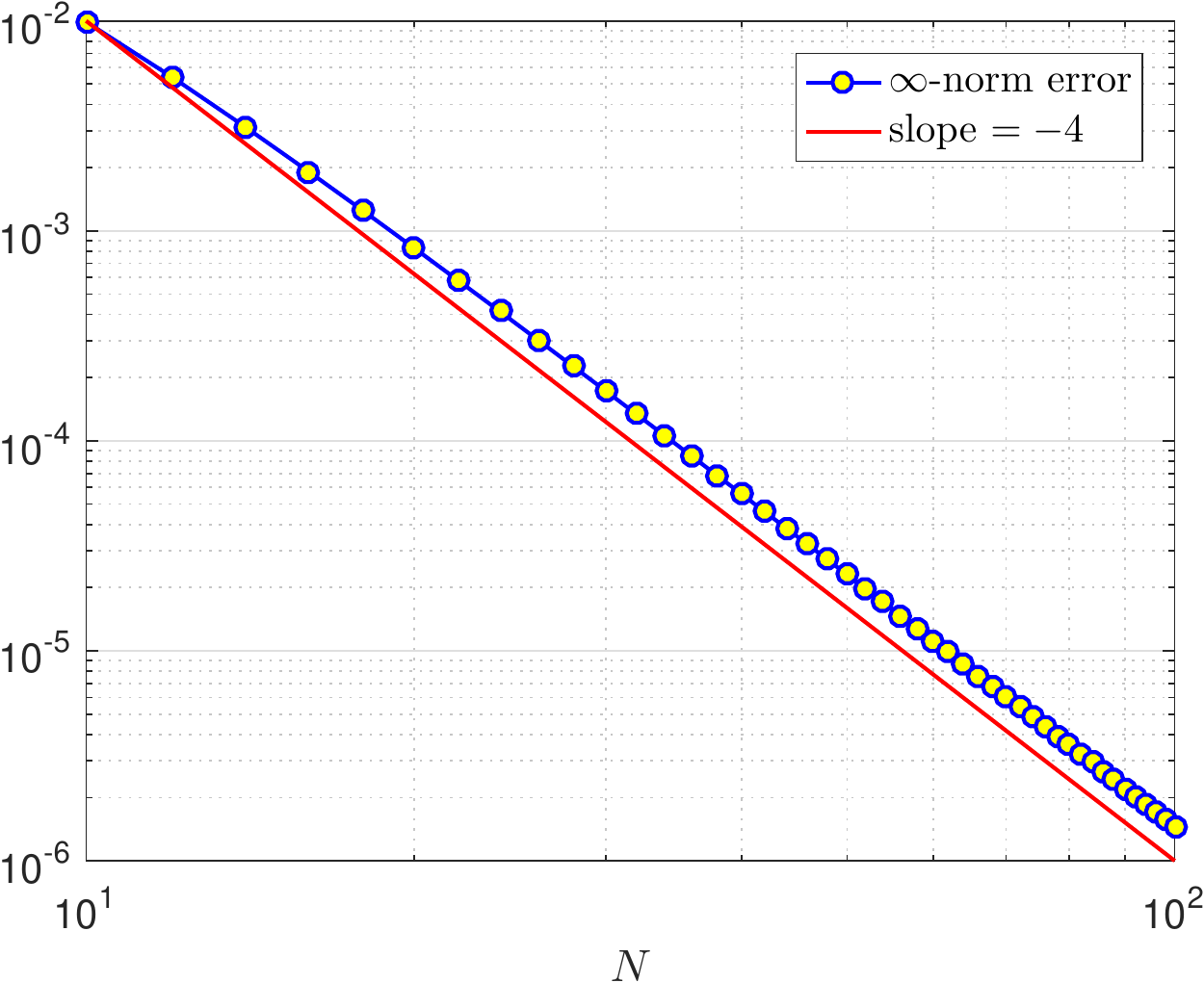}
  \caption{{\sc Example}~\ref{ex:abs}: Pointwise errors  $|f(x)-r(x)|$ of the  rational interpolant \eqref{eq:rat} for  $f(x)=|x|$ at \texttt{xx} with $N=40$ and $s=2$ (left),  and the convergence rates $\|f(\texttt{xx})-r(\texttt{xx})\|_{\infty}$ for $N=10:2:100$ (right).}
  \label{fig:ex1_abs}
\end{figure}

\begin{figure}[pt]
  \centering
  \includegraphics[width=.45\textwidth]{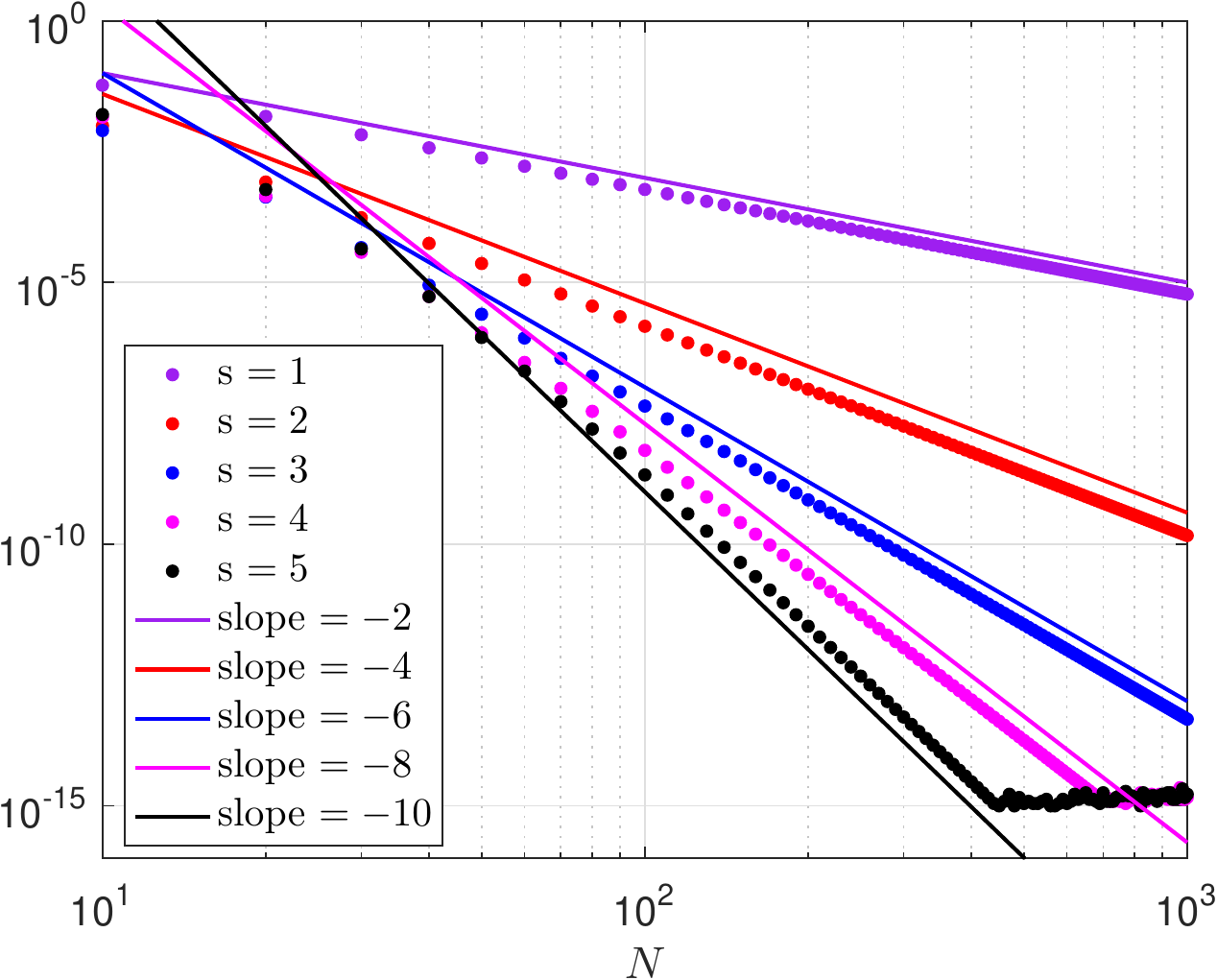}
  \includegraphics[width=.45\textwidth]{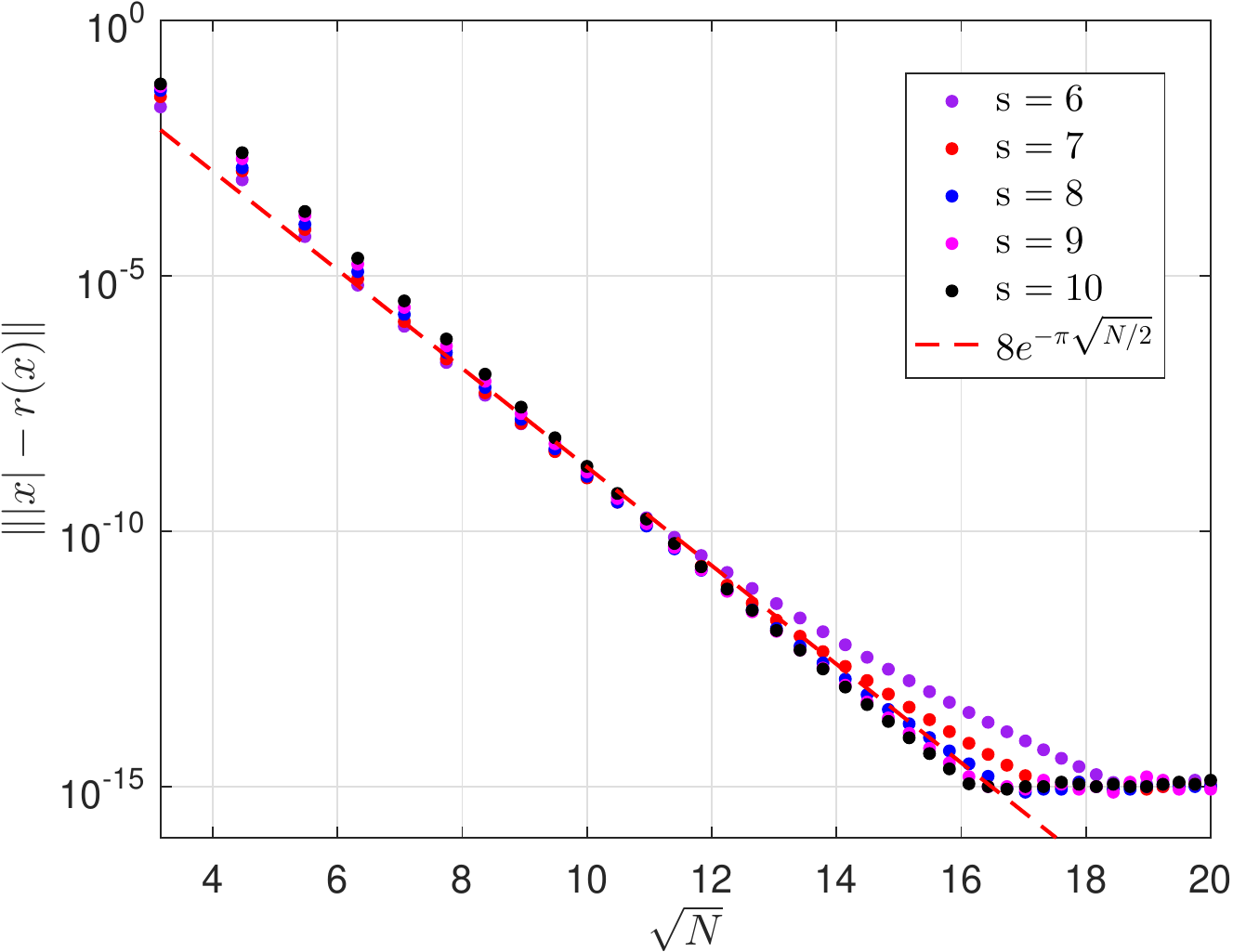}
  \caption{{\sc Example}~\ref{ex:abs}: Convergent rates  $\|f(\texttt{xx})-r(\texttt{xx})\|_{\infty}$ for  $f(x)=|x|$ with $N=10:10:1000$ for $s=1,2,\ldots,5$ (left) and $N=10:10:400$ for $s=6,7,\ldots,10$ (right). The left column of the figure is plotted in loglog axis, while the right one is shown in log scale with $\sqrt{N}$ on the horizontal axis. The asymptotically straight line on this axis shows the root-exponential effect. The red dashed line denotes the theoretical result by \eqref{eq:newmann2} with $\alpha=1$ and $N/2$.}
  \label{fig:ex1_abs_s}
\end{figure}

\begin{figure}[pt]
  \centering
  \includegraphics[width=.45\textwidth]{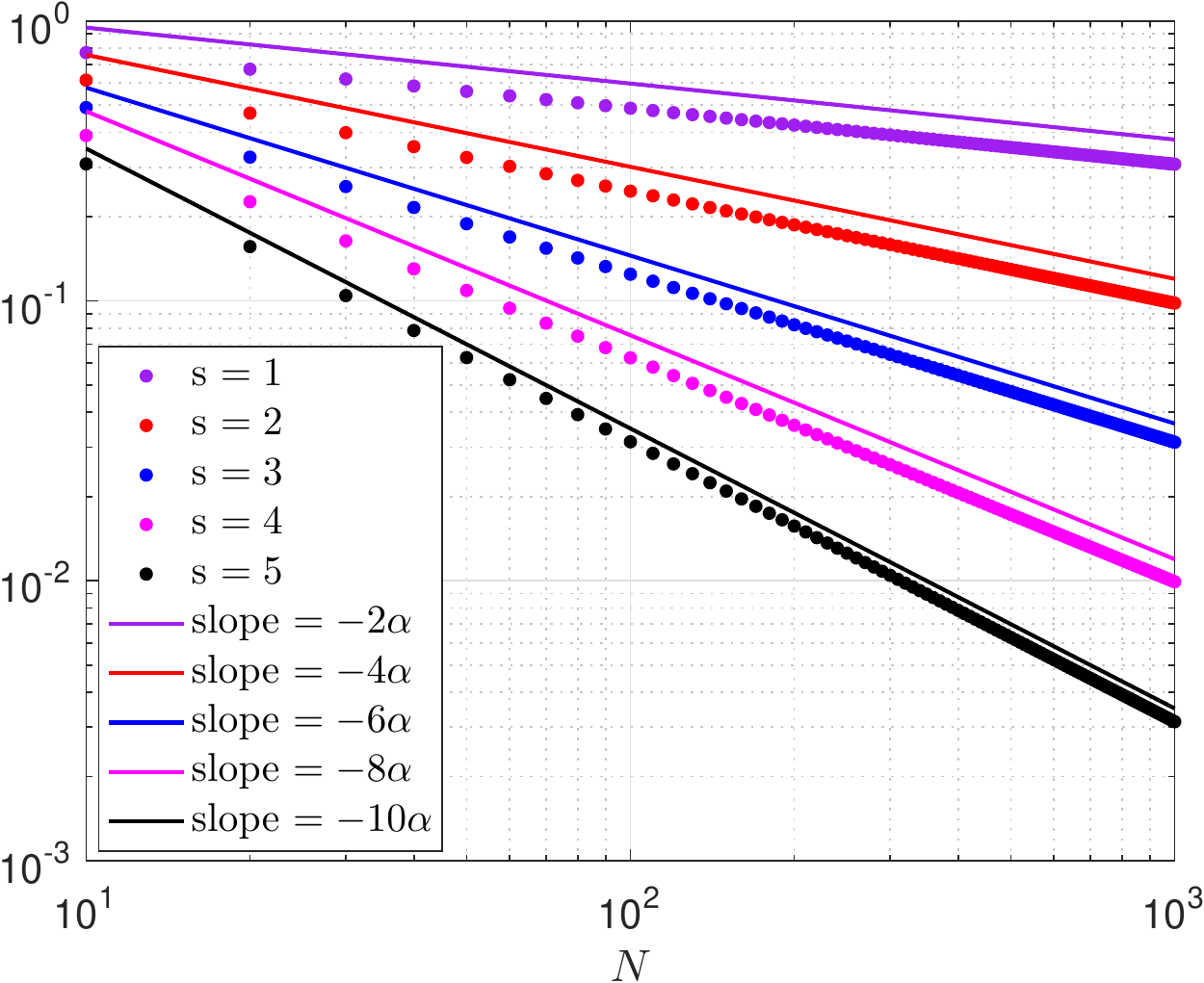}
  \includegraphics[width=.45\textwidth]{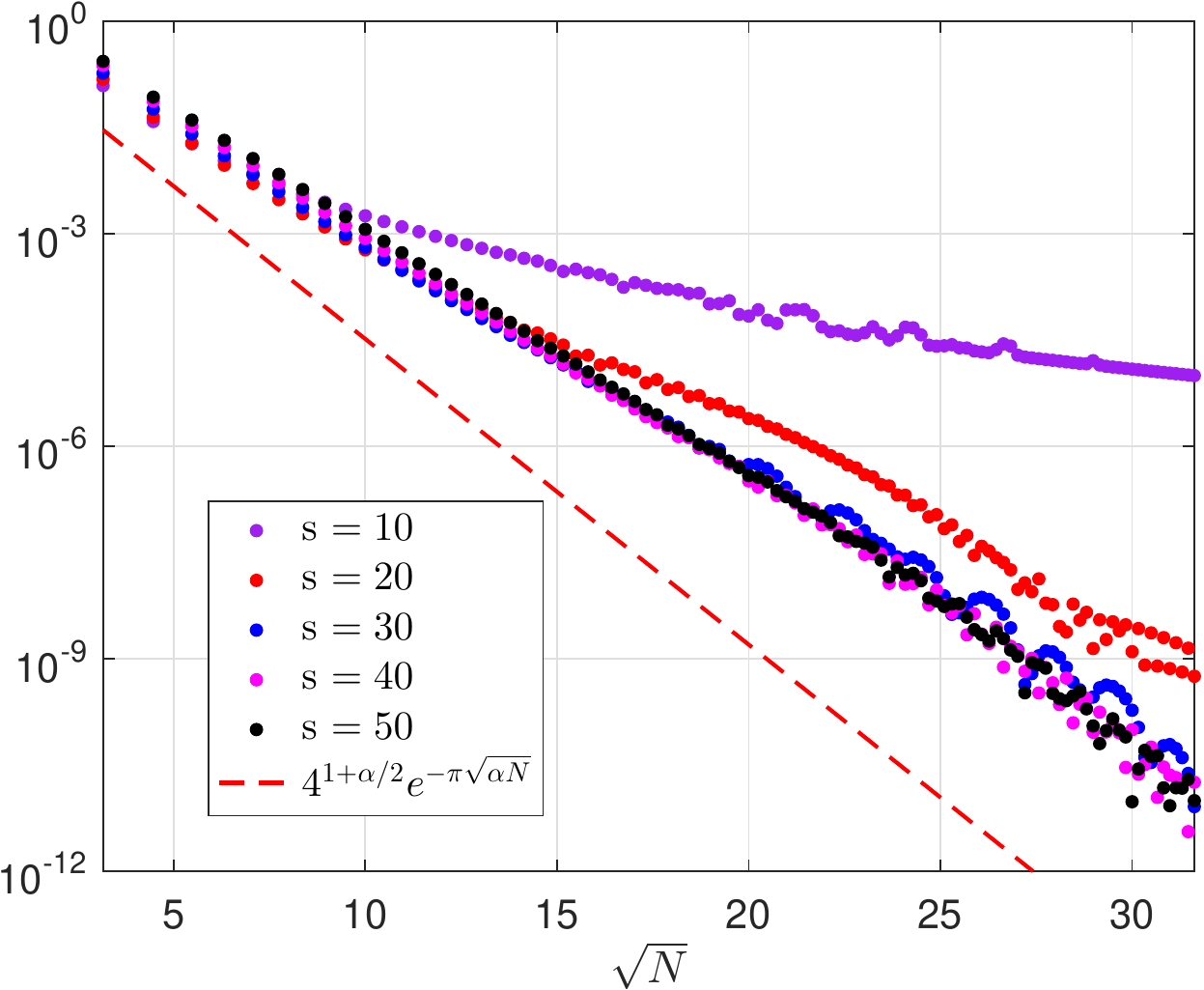}
  \caption{{\sc Example}~\ref{ex:abs}: Images as in {\sc fig.}~\ref{fig:ex1_abs_s} but now for function $|x|^\alpha$: rational interpolation to $|x|^{0.1}$ by taking $N=10:10:1000$ for $s=1,2,\ldots,5$ on the left and for $s=10,20,\ldots,50$ on the right.}
  \label{fig:ex1_abs_a_s}
\end{figure}

The left column of {\sc Fig.}~\ref{fig:ex1_abs} shows the pointwise error of $f(x)-r(x)$ at nodes \texttt{xx}. It is plain that the linear rational interpolant \eqref{eq:rat} performs well with a small error at the nonsmooth point $x=0$ and interpolates better as $|x| \to 1$. The right column of {\sc Fig.}~\ref{fig:ex1_abs} shows the errors in the infinite norm against $N$,
and illustrates the error decays  at a rate $N^{-4}$. 

In addition, the experimental results on the left column of {\sc Fig.}~\ref{fig:ex1_abs_s} show that the rate is algebraic $N^{-2s}$ for $s$ no more than 5.
As $s$  increases from  $8$ to $10$, the right column of {\sc Fig.}~\ref{fig:ex1_abs_s} illustrates that  the root-exponential rate is recovered.  In particular,  the approximation  \eqref{eq:rat}  converges as the theoretically best convergence rate \eqref{eq:newmann2} with $N/2$ instead of $N$.


%

Furthermore, an extensional experiment is on the rational approximation for $f(x)=|x|^\alpha$ for $\alpha \in (0,1)$ and $x\in [-1,1]$, whose theoretically best uniform  convergence rate is  given by \eqref{eq:newmann2}. To test the approximation property of the proposed rational interpolation, we consider  $\alpha \in (0,1)$ with  $\alpha=0.1$ (see {\sc Fig.}~\ref{fig:ex1_abs_a_s}).

{\sc Fig.}~\ref{fig:ex1_abs_a_s} shows the convergence rates, from which we know that an algebraical decay is obtaied with order $N^{-2\alpha s}$ for some values of $s$ and $N$. As $s$ increases from $10$ to $50$ in \eqref{eq:refpts2} (corresponding to $s$ from $1$ to $5$ in \eqref{eq:refpts}),  a root-exponential rate is recovered again. A phenomenon demonstrated by {\sc Fig.}~\ref{fig:ex1_abs_a_s} is that as $s$ increases, the root-exponential rate is achieved asymptotically from the algebraic rate. However, the transition between them is not clear.

From the equivalence between $|x|^\alpha$ defined on $[-1,1]$ and $x^\alpha$ defined on $[0,1]$ \cite{Stahl2003},  similar results can  be obtained for the rational interpolants for approximation of $x^\alpha$.



\begin{remark}
Similar to  the linghtning method  \cite{Gopal2019,Trefethen2020}  with  the poles  clustering exponentially near the singularity,
a lot of poles of the proposed rational interpolant also cluster near the origin (see {\sc Fig.}~\ref{fig:pol_zer}).
%
%
\begin{figure}[hptb]
  \centering
  \includegraphics[height=4cm,width=6cm]{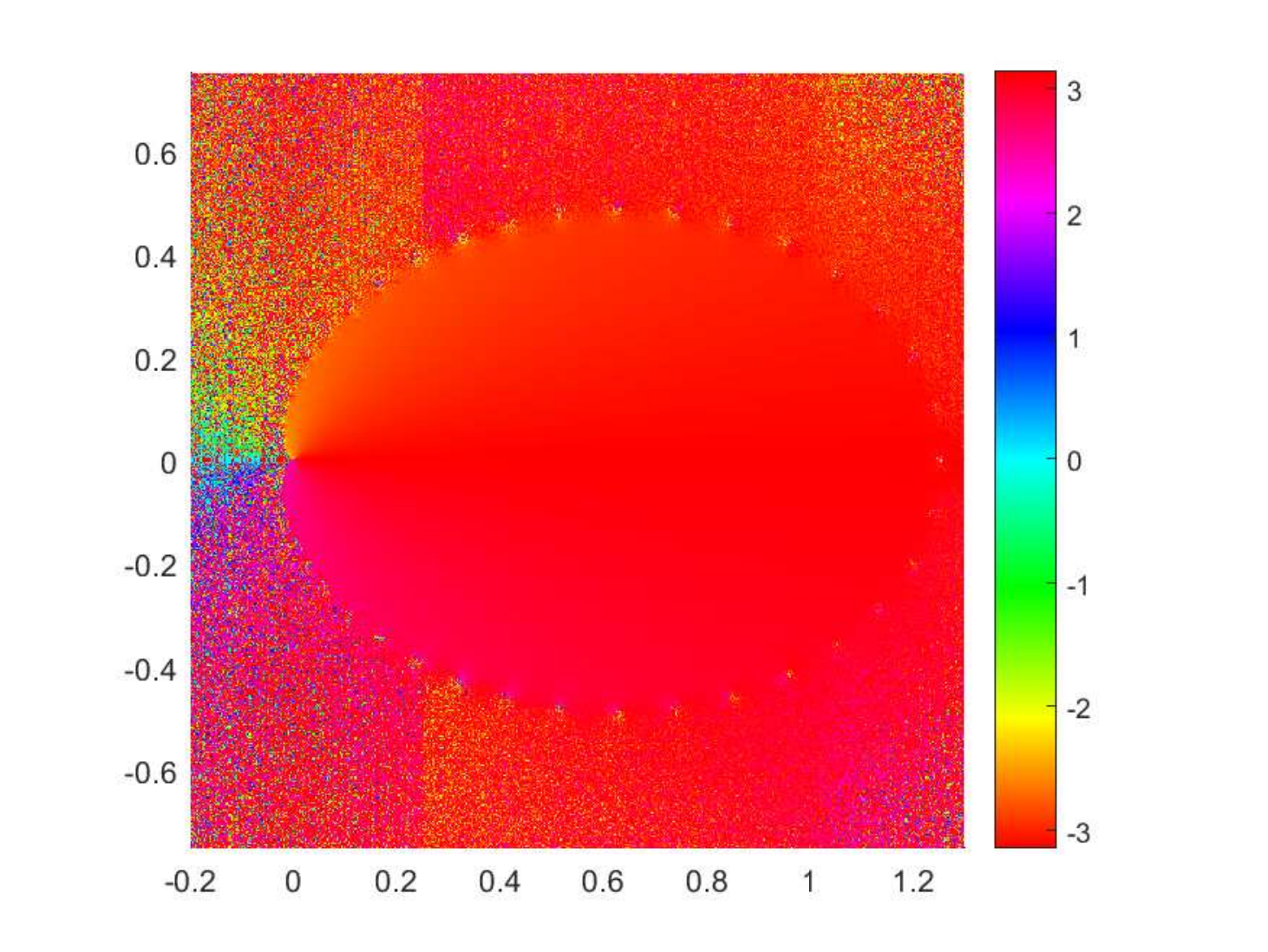}
  \includegraphics[height=3.85cm,width=4cm]{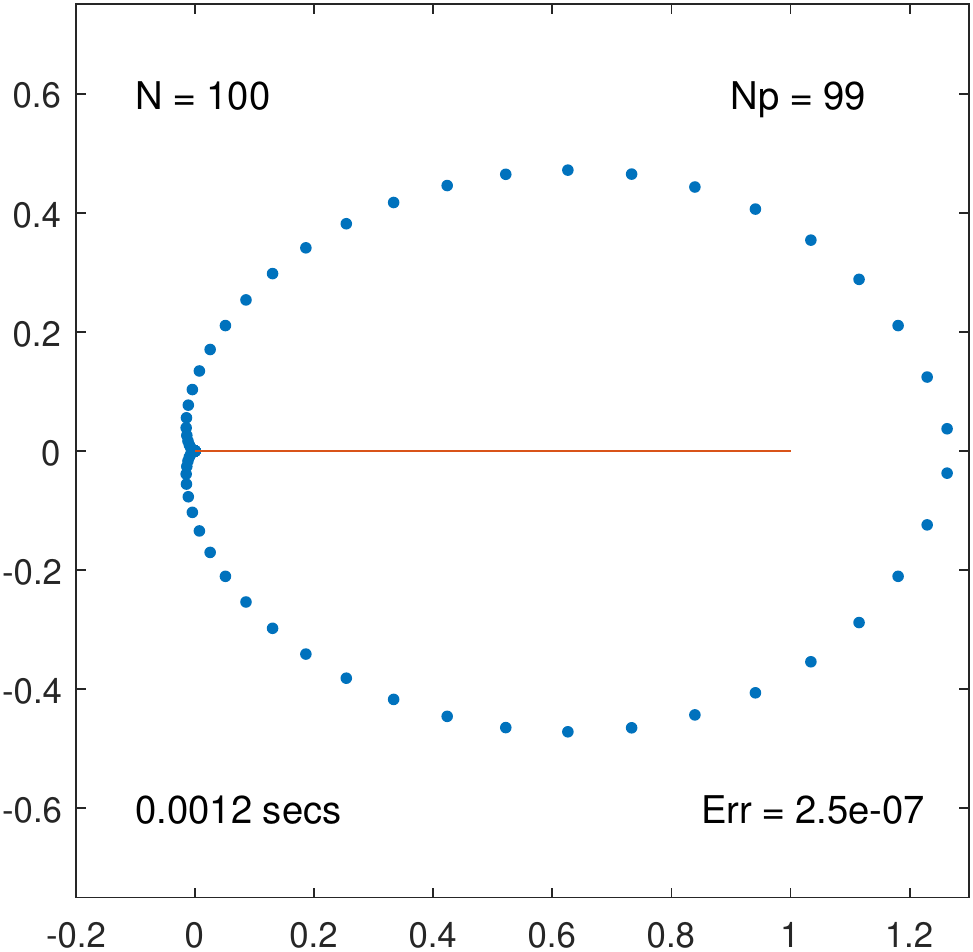} \\
  \includegraphics[height=4cm,width=6cm]{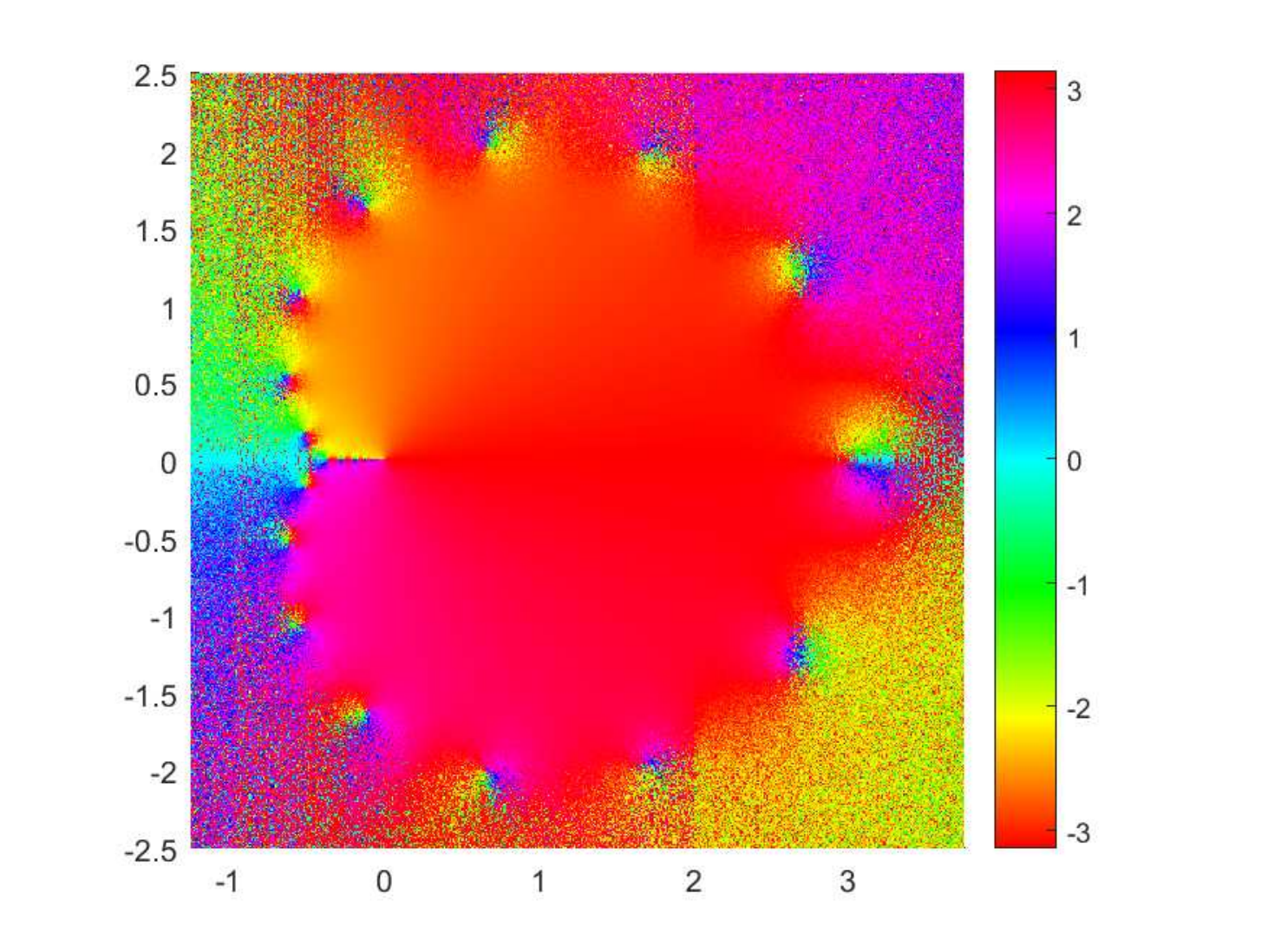}
  \includegraphics[height=3.85cm,width=4cm]{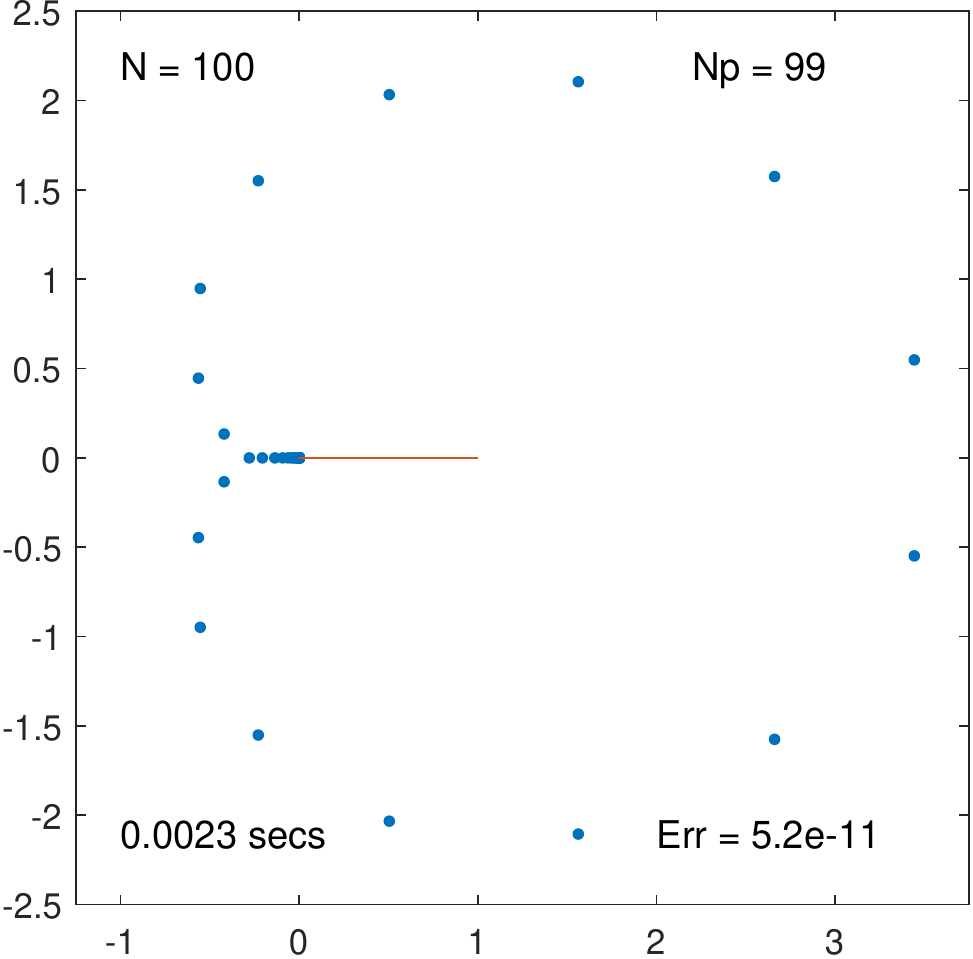} \\
  \includegraphics[height=4cm,width=6cm]{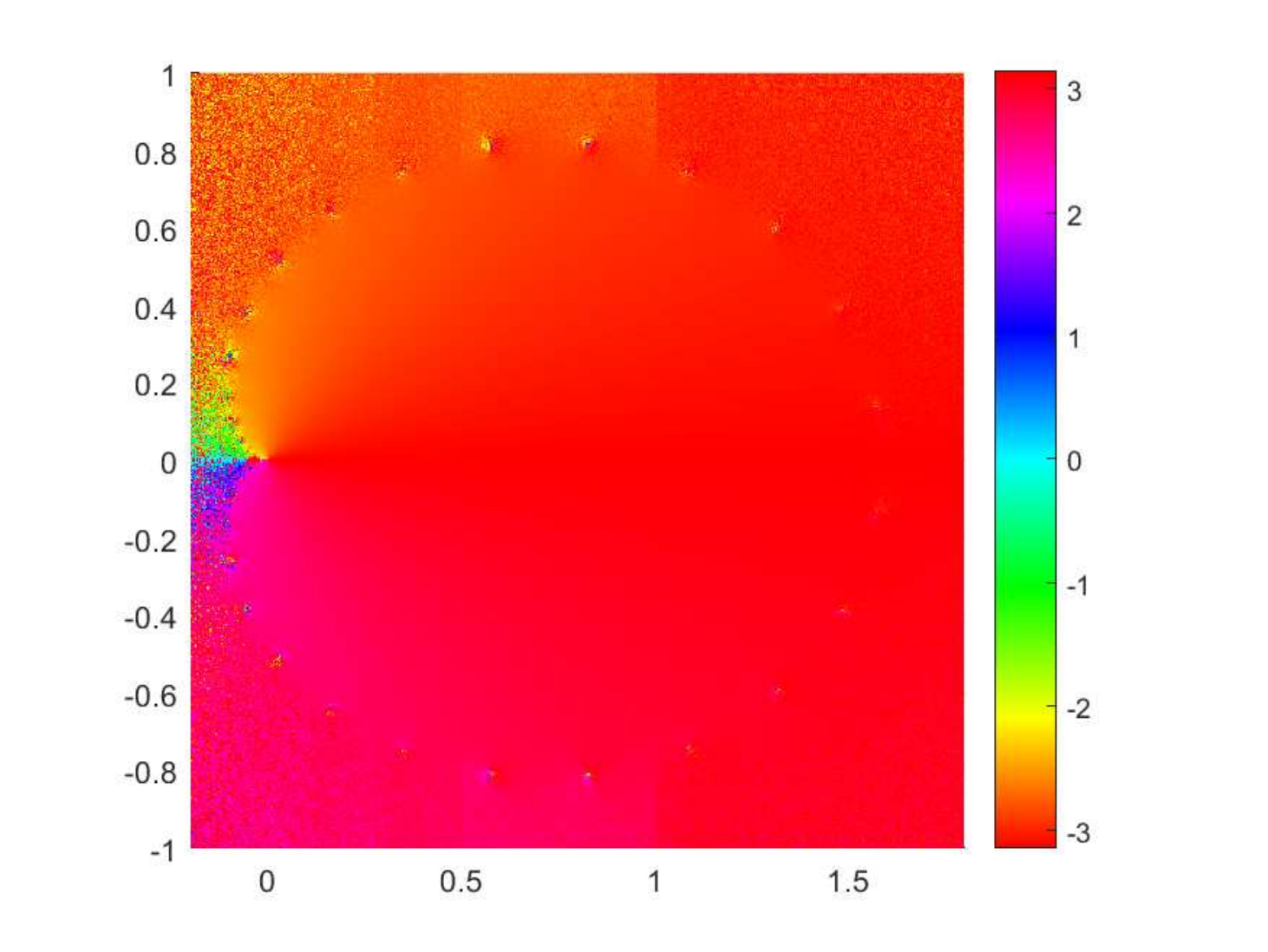}
  \includegraphics[height=3.85cm,width=4cm]{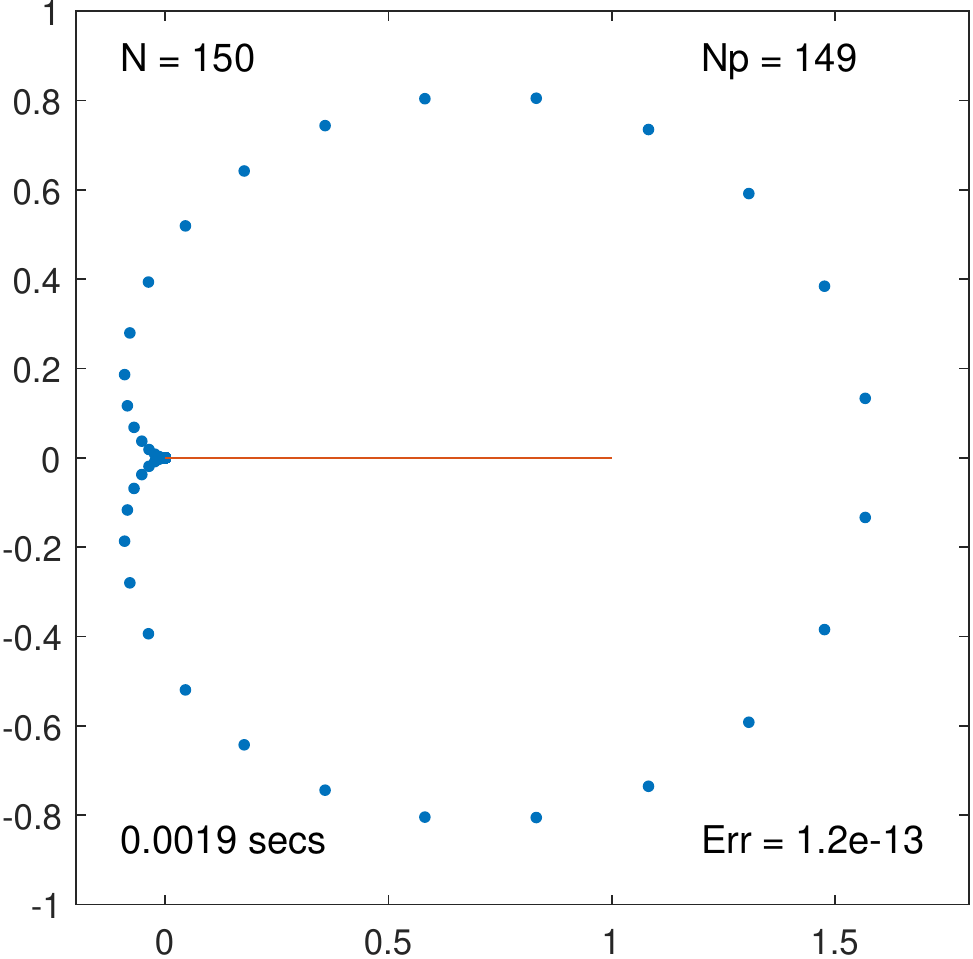}
  \caption{{\sc Example}~\ref{ex:abs}: Phase portraits (left column) and distributions of poles (right column) of rational interpolant \eqref{eq:rat} to $x^{1/\pi}$ on $[0,1]$ (red line) with $s = 2$ (Top) and $s = 10$ (Middle, Bottom). The number of interpolation nodes, \texttt{N}, and poles, \texttt{Np}, are listed on the upper left and upper right, respectively. The poles are numerically computed by Chebfun function \texttt{prz}. The lower left lists the elapsed time for computing the approximation. The number \texttt{Err} is equal to the maximum of $|f(x)-r(x)|$ over the discrete grid \texttt{xx}.}
  \label{fig:pol_zer}
\end{figure}
\end{remark}
\end{example}

\begin{example}[A further comparison of the three kinds of approximations for $f(x)=\sqrt{x}$: Newman's method \cite{Newman1964}, AAA method \cite{Nakat2018} and \eqref{eq:rat}]\label{ex:three}\rm
The Newman approximation is
\begin{align*}
  r(x)=\sqrt{x}\frac{p(\sqrt{x})-p(-\sqrt{x})}{p(\sqrt{x})+p(-\sqrt{x})},
\end{align*}
where $p(x)=\prod_{k=0}^{2N-1} (x+\xi^k)$ and $\xi=\exp(-(2N)^{-1/2})$. We apply AAA method with a set of 10,000 Chebyshev points in $[0,1]$. All errors are taken as the infinity norm in discrete set \texttt{xx}.

Results shown in {\sc Fig.}~\ref{fig:ex1_abs_comp} illustrate that the rational interpolant \eqref{eq:rat} behaves better than Newman approximation and the AAA algorithm as $N$ becoming larger. 

\begin{figure}[hbpt]
  \centering
  \includegraphics[width=.425\textwidth]{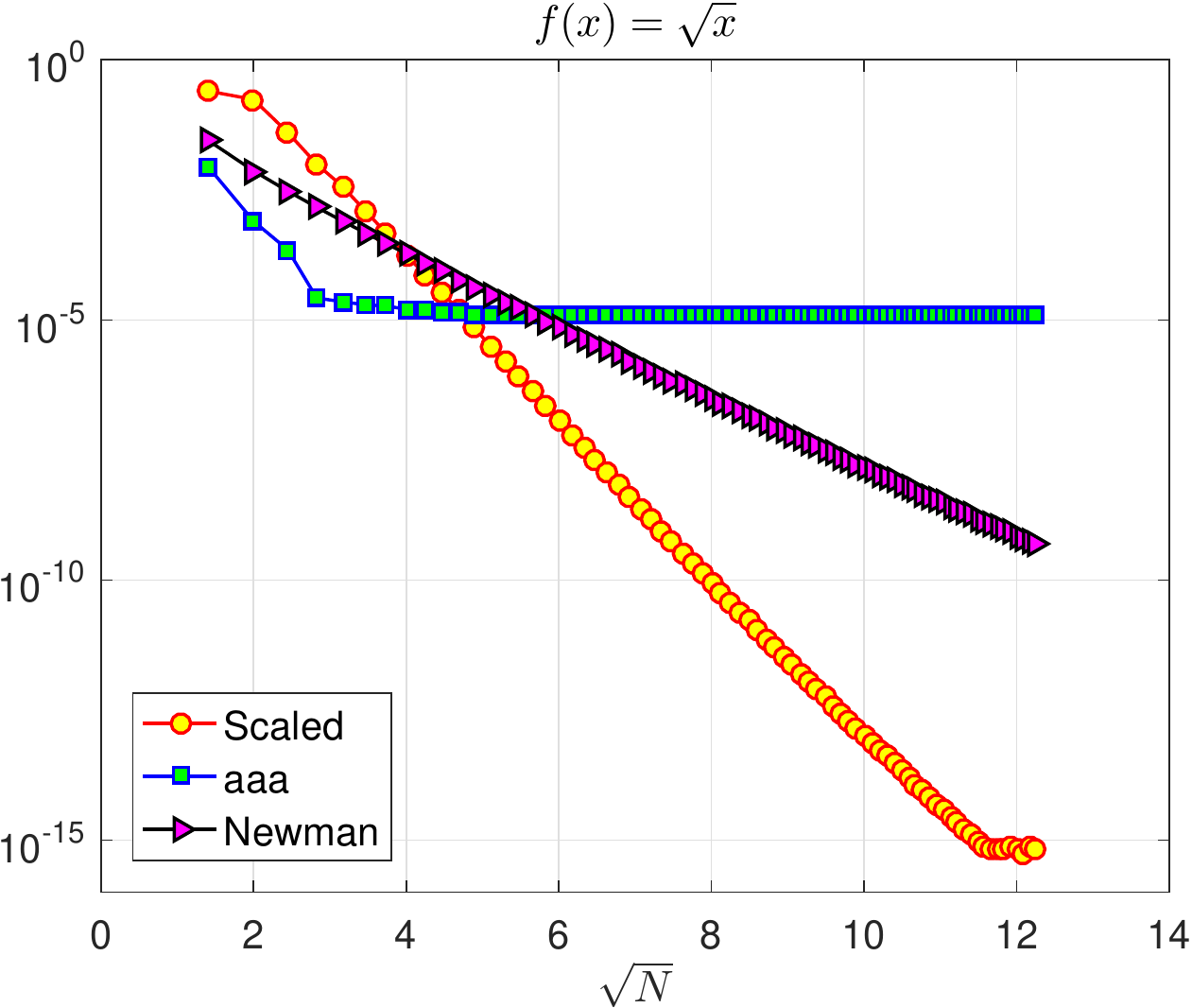}
  \caption{{\sc Example}~\ref{ex:three}: Comparison of three kinds of methods for $f(x)=\sqrt{x}$ on $[0,1]$: Newman's method, AAA method and \eqref{eq:rat} (denoted by `Scaled') with $s = 20$ for even values of $N$ from $2$ to $150$ in maximum error in \texttt{xx}. }
  \label{fig:ex1_abs_comp}
\end{figure}
\end{example}

\begin{example}[Volterra integral equation with a weakly singular kernel]\label{ex:vol}\rm
  In this example, we intend to study the numerical solution of Volterra integral equation of the second kind:
  \begin{align}\label{eq:vol}
    y(x) = f(x) + \int_0^x (x-t)^{-\alpha} K(x,t) y(t) dt, \quad x \in \Omega
  \end{align}
  for $\alpha \in (0,1)$, $f(x)$ and $K(x,t)$ defined respectively on $\Omega$ and $\mathcal{D}:=\{ (x,t) : 0 \leq t \leq x \leq T \}$, and $K(x,x) \neq 0$ for $x \in \Omega$.
  The regularity of the solution has been extensively investigated: equation \eqref{eq:vol}  has a unique solution $y \in C^m(0,T] \cap C(\Omega)$ with $|y'(x)| \leq Cx^{-\alpha}$ provided $f \in C^m(\Omega)$ and $K \in C^m(\mathcal{D})$ for some $m \geq 1$ \cite{Brunner2004}.

  Since the singularity of the solution $y(x)$ at the origin, the collocation method with piecewise polynomial on graded meshes
  is preferred \cite{Liang2019}. Here we utilize the linear rational interpolation to construct a global approximate solution.
  Let us firstly discretize the successive domain $\Omega$ into $N$ parts with nodes $\{ x_k \}_{k=0}^N$ given by \eqref{eq:refpts1}.
  The numerical scheme can be written as
  \begin{align*}
    y_N(x_k) &= f(x_k) + \sum_{i=0}^N \lambda_i \lrbrack{ \int_0^{x_k} \frac{(x_k-t)^{-\alpha} K(x_k,t)}{(t-x_i) \sum_{l=0}^N \lambda_l /(t-x_l)} dt } y_N(x_i) \\
    &= f(x_k) + \sum_{i=0}^N \lambda_i \lrbrack{ (x_k)^{1-\alpha} \int_0^{1} \frac{(1-\theta)^{-\alpha} K(x_k,x_k \theta)}{(x_k \theta-x_i) \sum_{l=0}^N \lambda_l /(x_k \theta-x_l)} d\theta } y_N(x_i) \\
    &= f(x_k) + \lambda_0 \lrbrack{ (x_k)^{1-\alpha} \int_0^{1} \frac{(1-\theta)^{-\alpha} K(x_k,x_k \theta)}{(x_k \theta-x_0) \sum_{l=0}^N \lambda_l /(x_k \theta-x_l)} d\theta } y_N(x_0) \\
    &\quad + \sum_{i=1}^N \lambda_i \lrbrack{ (x_k)^{1-\alpha} \int_0^{1} \frac{(1-\theta)^{-\alpha} K(x_k,x_k \theta)}{(x_k \theta-x_i) \sum_{l=0}^N \lambda_l /(x_k \theta-x_l)} d\theta } y_N(x_i),
  \end{align*}
which leads to  the matrix form
  \begin{align}\label{eq:volm}
    (\mathbf{I}-\mathbf{B})\mathbf{y} = \mathbf{b},
  \end{align}
  where $\mathbf{y} = (y_N(x_1), y_N(x_2), \ldots, y_N(x_N))^\top$, $\mathbf{b} = (b_1, b_2, \ldots, b_N)^\top$ with
  \begin{align*}
    b_i = f(x_i) + \lambda_0 \lrbrack{ (x_i)^{1-\alpha} \int_0^{1} \frac{(1-\theta)^{-\alpha} K(x_i,x_i \theta)}{(x_i \theta-x_0) \sum_{l=0}^N \lambda_l /(x_i \theta-x_l)} d\theta } y_N(x_0),\ i = 1,\ldots,N,
  \end{align*}
  $\mathbf{I}$ is  the identity matrix with dimension $N \times N$, and
  \begin{align*}
    \mathbf{B} = (B_{ki}) =
    \begin{pmatrix}
      \displaystyle (x_k)^{1-\alpha} \int_0^{1} \frac{(1-\theta)^{-\alpha} K(x_k,x_k \theta)}{(x_k \theta-x_i) \sum_{l=0}^N \lambda_l /(x_k \theta-x_l)} d\theta \\
      \\
      i,k = 1,\ldots,N
    \end{pmatrix}.
  \end{align*}
  In numerical implementation, we set $y_N(x_0) = f(x_0)$ and compute the entries of matrix $\mathbf{B}$ by numerical quadrature formula.
  Once \eqref{eq:volm} has been solved, the evaluation of $y(x)$ associated with a set of points $\{x_k\}_{k=0}^N$ amounts to substituting $\mathbf{y}$ into the rational interpolation formula \eqref{eq:rat}.

  Consider $f(x)=1$, $K(x,t)=\frac{1}{10\Gamma(1-\alpha)}$ such that the exact solution can be explicitly expressed as $y(x) = E_{1-\alpha,1}(x^{1-\alpha}/10)$, where the Mittag-Leffler function $E_{\mu,\nu}$ is defined by
  \begin{align*}
    E_{\mu,\nu}(z) := \sum_{p=0}^\infty \frac{z^p}{\Gamma(\mu p + \nu)} \quad \text{ for } \quad \mu,\nu,z \in \mathbb{R} \quad \text{ with } \quad \mu>0.
  \end{align*}
  Numerical results shown in {\sc Fig.}~\ref{fig:vol} for different values of $\alpha$ and $s$ illustrate that a root-exponential convergence rate is achieved by varying $N$, the number of interpolation nodes.

  \begin{figure}[tp]
    \centering
    \includegraphics[width=.6\textwidth]{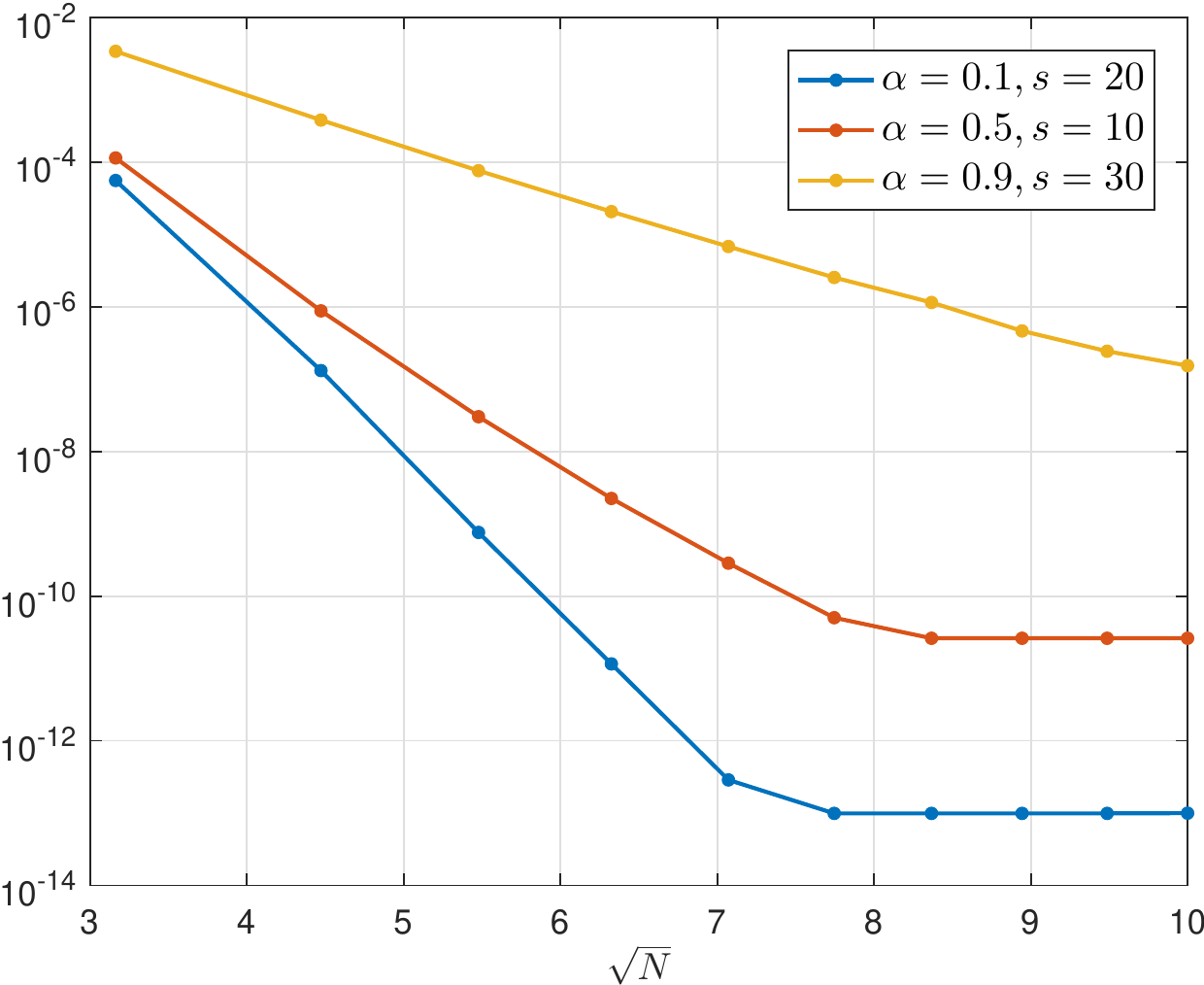}
    \caption{{\sc Example}~\ref{ex:vol}: Approximation errors in the infinite norm at \texttt{xx} of \eqref{eq:rat} to problem \eqref{eq:vol} with $N=10:10:100$ and different values of $\alpha$ and $s$.}
    \label{fig:vol}
  \end{figure}
\end{example}

\vspace{0.4cm}
\begin{remark}
  The linear barycentric rational interpolation \eqref{eq:rat} is also available for the situations where Jacobi-Gauss-Lobatto nodes and simplified weights are used.
Next, we shall select $\{ \lambda_i \}_{i=0}^N$ as the simplified weights corresponding to the Jacobi-Gauss-Lobatto nodes for $\{ y_i \}_{i=0}^N$, and determine $\{ x_i \}_{i=0}^N$ by mapping $\{ y_i \}_{i=0}^N$ with respect to \eqref{eq:refpts}.

Note that the simplified weights associated with Jacobi-Gauss-Lobatto points, zeros of $(1-x^2) \partial_x J_N^{(\beta, \gamma)}(x)$, are proposed in \cite[Corollary 3.3]{Wang2014}:
\begin{align}\label{eq:jacweig}
  \lambda_i = (-1)^i \sqrt{\delta_i \omega_i},\quad \text{for} \quad
  \delta_i =
  \begin{dcases}
    1, & i=1,2\ldots,N-1, \\
    \gamma+1, & i=0, \\
    \beta+1, & i = N.
  \end{dcases}
\end{align}
In \eqref{eq:jacweig}, $\{ \omega_i \}_{i=0}^N$ are the corresponding weights of the interpolatory quadrature rule corresponding to the weight function $(1-x)^\beta (1+x)^\gamma$.
With this choice, in the following example, we will illustrate numerically that the proposed linear barycentric rational formula \eqref{eq:rat} is experimentally well-conditioned, too.

\begin{example}[Highly oscillatory problem]\label{ex:hop}\rm
  One application considered in this example is about solving a highly oscillatory problem with a singulary kernel \cite{Wang2020a}:
  \begin{align}\label{eq:hop}
    I_\omega^{[0,a]}[f,h] := \int_0^a f(x) x^\alpha e^{i\omega h(x)} dx, \quad \alpha \in (0,1) \mbox{ and } \omega \gg 0,
  \end{align}
  via Levin's method. This method reads: Find a function $p$ such that
  \begin{align*}
    p'(x) e^{i\omega h(x)} + i\omega h'(x) p(x) e^{i\omega h(x)} = (p(x) e^{i\omega h(x)})' = f(x) x^\alpha e^{i\omega g(x)},
  \end{align*}
  which is equivalent to
  \begin{align}\label{eq:leven}
    p'(x) + i\omega h'(x) p(x) = f(x) x^\alpha,
  \end{align}
  then, the highly oscillatory problem can be computed by
  \begin{align*}
    I_\omega^{[0,a]}[f,h] = p(a) e^{i\omega h(a)} - p(0) e^{i\omega h(0)}.
  \end{align*}
  In \eqref{eq:hop}, we assume functions $f$ and $h$ are both smooth and $h'(x) \neq 0$ on $\Omega$. The main advantages of Levin's method is converting equivalently the integral with a highly oscillatory kernel to a non-oscillatory ordinary differential equation (ODE) problem whose approximation determines the accuracy of \eqref{eq:hop} \cite{Levin1996}.  In \eqref{eq:leven} for $\alpha \in (0,1)$, $p(x)$ has singularity at $x=0$ because of the singular source term $f(x)x^\alpha$. Thus, the appropriate method should be carefully designed to catch the singularity.

  Applying the rational approximation proposed in this paper to \eqref{eq:leven} with the differential matrix $D^{(1)}$ given by \eqref{eq:diffm1}, we obtain a scheme by collocation method written in matrix form as
  \begin{align}\label{eq:coll}
    (D^{(1)} + i\omega h'(\mathbf{x})) \mathbf{p} = \mathbf{f}
  \end{align}
  with
  \begin{gather*}
    \mathbf{x} = (x_0, x_1, \ldots, x_N)^\top, \quad
    \mathbf{p} = (p(x_0), p(x_1), \ldots, p(x_N))^\top, \\
    \mathbf{f} = (f(x_0)x_0^\alpha, f(x_1)x_1^\alpha, \ldots, f(x_N)x_N^\alpha)^\top.
  \end{gather*}
  Once \eqref{eq:coll} has been solved, the evaluation of $I_\omega^{[0,a]}[f,h]$ amounts to the computation of the following formula
  \begin{align}
    I_{N,\omega}^{[0,a]}[f,h] := \mathbf{p}(N+1) e^{i\omega h(a)} - \mathbf{p}(1) e^{i\omega h(0)},
  \end{align}
  where $\mathbf{p}(i)$ denotes the $i$-th component of vector $\mathbf{p}$.

  Setting $f(x) = e^{i\omega}(1-x)(2-x)^\alpha$, $h(x)=x$ and $a=1$ in \eqref{eq:hop},
  we can obtain a reference solution by calling MATLAB code \texttt{quadgk} equipped with setting both relative error tolerance and absolute error tolerance to $10^{-14}$.
  In \cite{Wang2020a}, to avoid the singularity, the authors apply a singularity separation technique converting the singular ODE \eqref{eq:leven} into two kinds of non-singular ODEs.
Here we will approximate the derived ODE directly with rational interpolation \eqref{eq:rat}.
 As pointed before, the interpolation points can be chosen as  the Jacobi-Gauss-Lobatto nodes in $[-1,1]$ together with scaled map  \eqref{eq:refpts} or \eqref{eq:refpts2}. The weights in \eqref{eq:rat} are the simplified barycentric weights \eqref{eq:jacweig} corresponding to Jacobi-Gauss-Lobatto nodes in $[-1,1]$. In MATLAB, these points and weights can be implemented by calling \texttt{jacglquad}\footnote{\texttt{jagslb} is a MATLAB code computing the Jacobi-Gauss-Lobatto points on $[-1,1]$, and this code is available online: \url{https://www.ntu.edu.sg/home/lilian/book.htm}}:
\begin{verbatim}
  function [x, uk] = jacglquad(N, bet, gam, dom)
  % JACGLQUAD computes the Jacobi-Gauss-Lobatto points x in dom and
  % correspondingly simplified weights uk by calling CHEBFUN function
  % baryWeights.
  if bet == -0.5 && gam == -0.5
     [x,~,uk] = chebpts(N, [0,1]); x = max(dom)*x; return;
    else
     x = jagslb(N, bet, gam); x = max(dom)*(x+1)/2;
     uk = baryWeights(x);
  end
\end{verbatim}

  The simulated results shown in {\sc Fig.}~\ref{fig:levin} illustrate that the collocation method with rational interpolant is valid for fixed value $\omega$, and even converges in exponential rate against $N$.
  It should be noted that for fixed Jacobi-Gauss-Lobatto points, as  $s$ enlarges, the condition number of differential matrix \eqref{eq:diffm1} becomes large, since there are more collocation points accumulating at $x=0$ (singular point) that leads the reciprocal of $x_i-x_j$ getting larger for adjacent points $x_i$, $x_j$.

\begin{figure}[tp]
  \centering
  \includegraphics[width=.49\textwidth]{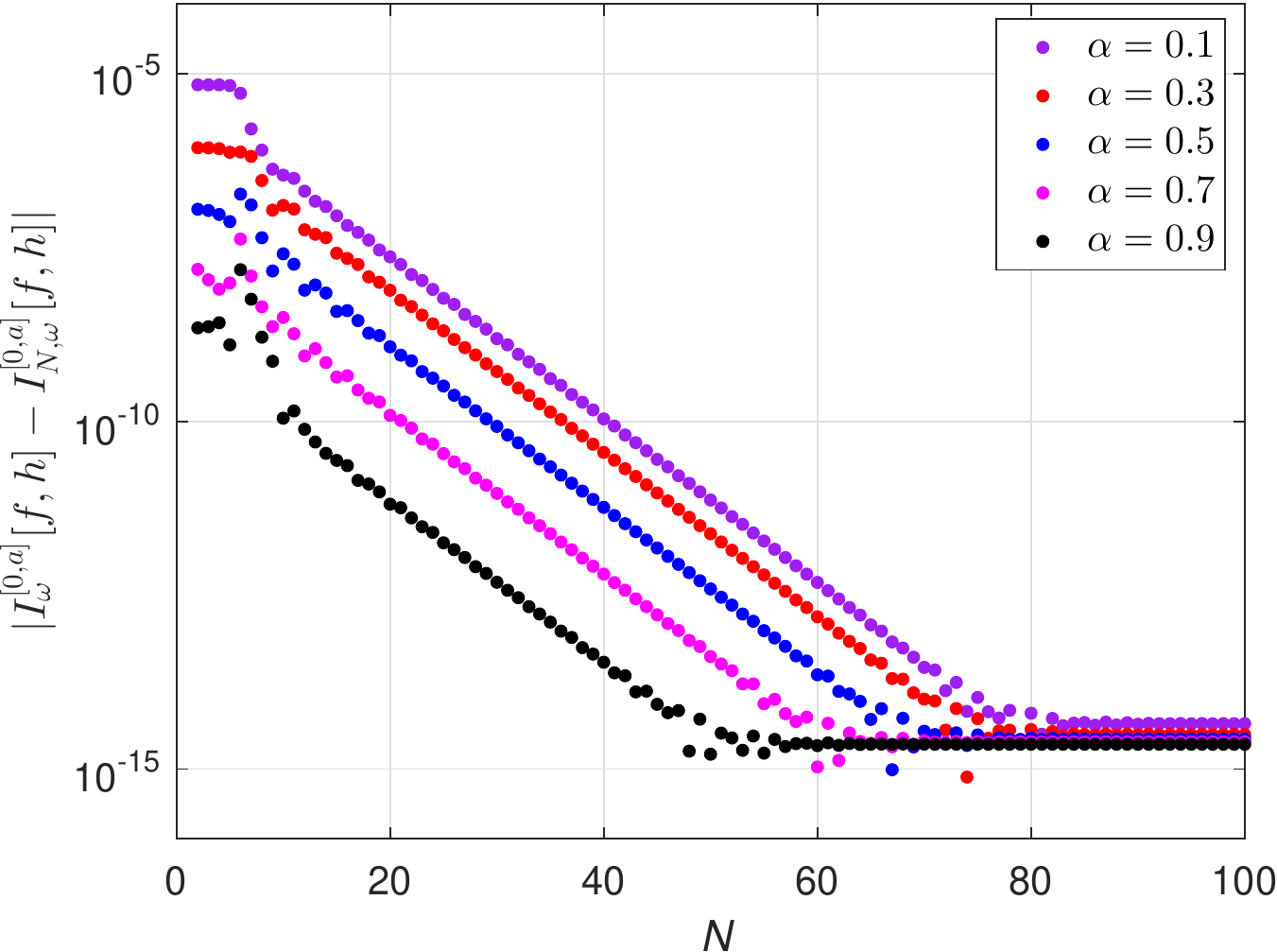}
  \includegraphics[width=.49\textwidth]{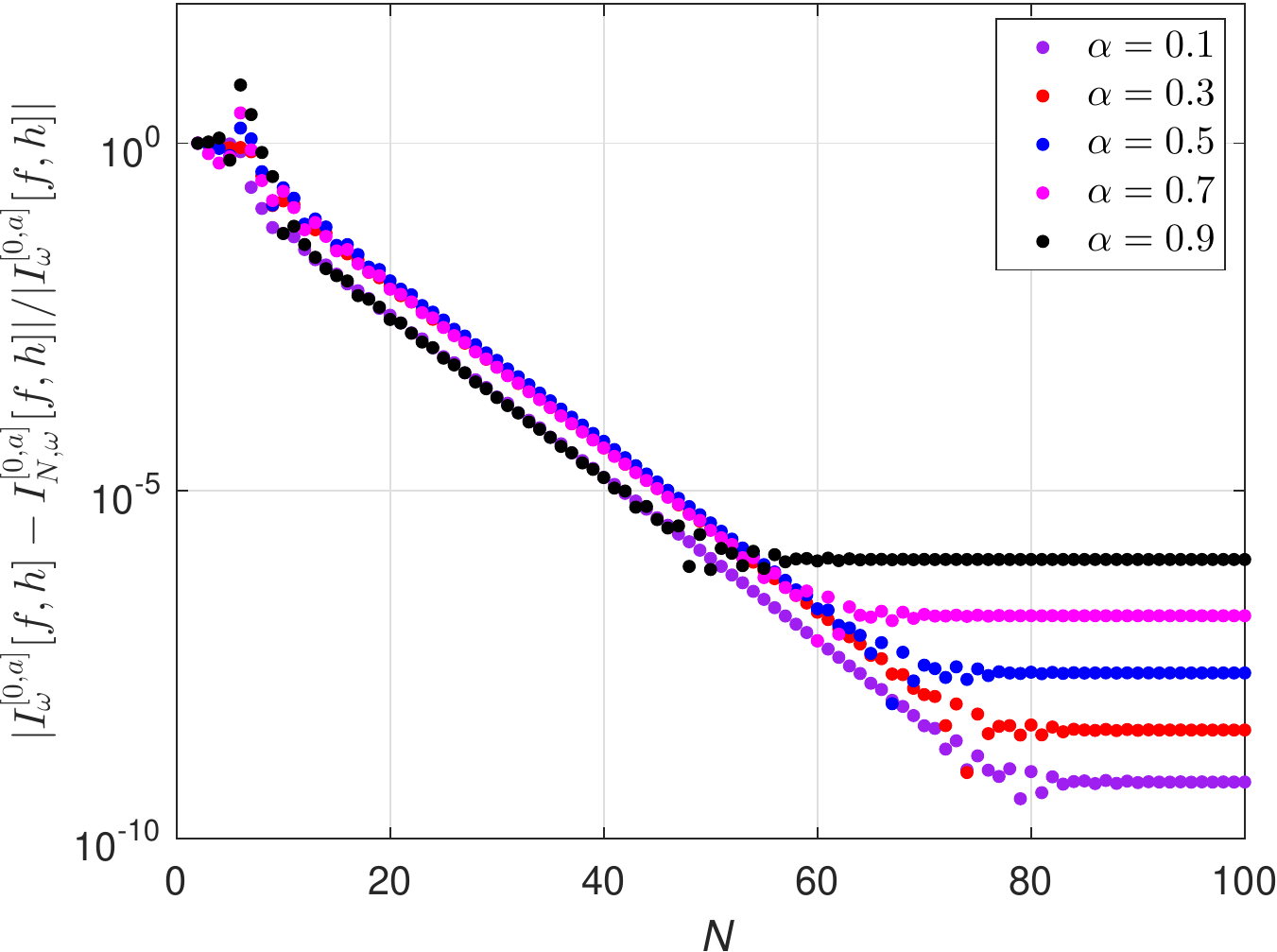}
  \caption{{\sc Example}~\ref{ex:hop}: Exponential convergence of the collocation method with  the  rational approximation against $N=2:100$. On the left column of the figure, the absolute errors with $\alpha$ taking $0,1, 0.3, \ldots, 0.9$ against $N$ are shown, while the right column of the images shows the relative errors. In this simulation, we take $s=10, \omega = 50,000$, interpolation nodes being Jacobi-Gauss-Lobatto points associated with Jacobi polynomial $J_N^{(0,10)}(y)$. }
   \label{fig:levin}
\end{figure}
\end{example}
\end{remark}


\section{Extension to functions of logarithmic singularity}\label{sec:exten}

\begin{figure}[hpt]
  \centering
  \includegraphics[width=.6\textwidth]{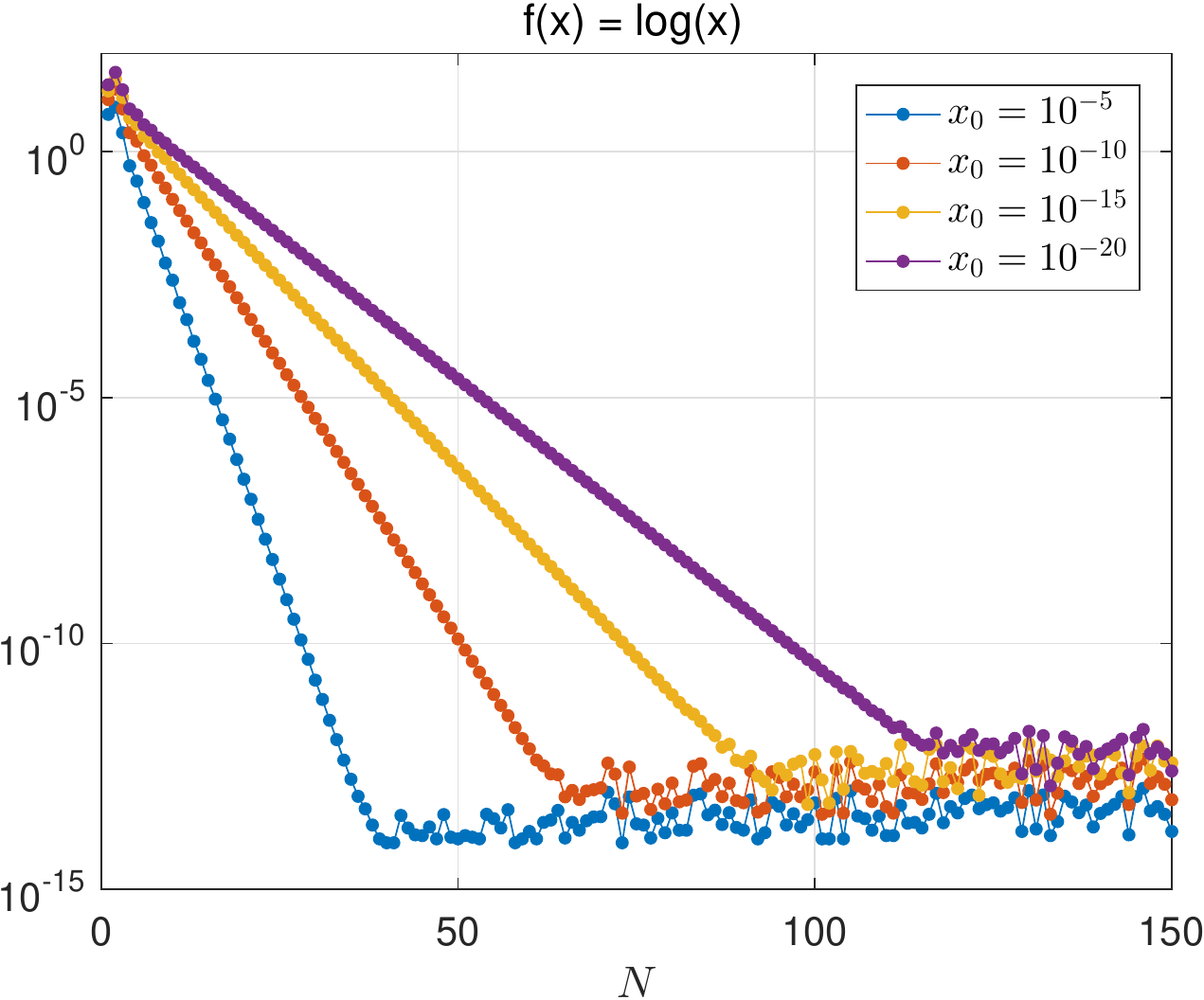}
  \caption{Rational approximation \eqref{eq:rat} with map \eqref{eq:refpts_log} to $\log(x)$ defined on $[x_0,1]$ with $x_0=10^{-20}, 10^{-15}, 10^{-10}$ and $10^{-5}$, respectively, against $N$ varying from $1$ to $150$. Maximum error takes the values in \texttt{xxx}.}
  \label{fig:ex:ratlog}
\end{figure}

In this section, we shall extend the rational approximation to functions of logarithmic singularity, i.e., $f(x) \sim \log(x)$ as $x \to 0^+$. Let us consider function $f(x)=\log(x)$ defined on domain $[x_0, T] \subset \Omega$ for some $x_0 > 0$.
The map defined in \eqref{eq:refpts_log} is equivalent to
\begin{align*}
  x = x_0 \lrbrack{ \frac{T}{x_0} }^{\frac{\theta+1}{2}} \quad \mbox{for} \quad \theta \in [-1,1].
\end{align*}
Calling \texttt{ratlog} in MATLAB (see {\sc Fig.}~\ref{fig:ratlog}) and taking maximum error in set \texttt{xxx}:
\begin{verbatim}
  xxx = logspace(log10(xmin), log10(xmax), 10000);
\end{verbatim}

We derive the rational approximation for functions of logarithmic singularity shown in {\sc Fig.}~\ref{fig:ex:ratlog} on approximation of $\log(x)$ for example. Numerical results illustrate the rational interpolant \eqref{eq:rat} with map \eqref{eq:refpts_log} to $\log(x)$ can achieve exponential convergence.

%

\section{Analysis of the rational approximation near the singularity}\label{sec:conv}
Now we start to analyse the accuracy  of the proposed rational interpolant. The basic tool for estimating the accuracy of rational approximation is the Hermite integral formula given in the following theorem.
\begin{theorem}[Hermite integral formula for rational interpolation \cite{Trefethen2020}]
  Let $\Omega$ be a simply connected domain in $\mathbb{C}$ bounded by a closed curve $\Gamma,$ and let $f$ be analytic in $\Omega$ and extend continuously to the boundary. Let interpolation points $\alpha_{0}, \ldots, \alpha_{N} \in \Omega$ and poles $\beta_{1}, \ldots, \beta_{N}$ anywhere in the complex plane be given. Let $r$ be the unique type $(N, N)$ rational function with simple poles at $\{\beta_{j}\}$ that interpolates $f$ at $\left\{\alpha_{j}\right\}$. Then for any $z \in \Omega$
  \begin{align*}
    f(z)-r(z)=\frac{1}{2 \pi i} \int_{\Gamma} \frac{\phi(z)}{\phi(t)} \frac{f(t)}{t-z} d t,
  \end{align*}
  where
  \begin{align*}
    \phi(z)=\left. \prod_{j=0}^{N} \lrbrack{ z-\alpha_{j} } \right/ \prod_{j=1}^{N} \lrbrack{ z-\beta_{j} }.
  \end{align*}
\end{theorem}
The use of this formula is based on the fact that one knows the locations and properties of the poles in advance \cite{Gopal2019,Trefethen2020} or approximate it accurately. For linear barycentric rational interpolant \eqref{eq:rat}, its poles can be approximated by solving the eigenvalues of a generalized companion matrix pair \cite{Lawrence2013}.

Another strategy in rational approximation is proposed in \cite[Theorem 4]{Balten1999}, where  the exponentially convergent rate for transformed Chebyshev points for analytic functions is established. A simplified version is stated in the following theorem.
\begin{theorem}[{\cite{Balten1999}}]\label{thm:cite1999}
  Let $\mathcal{D}_{1}, \mathcal{D}_{2}$ be two domains of $\mathbb{C}$ containing $J=[-1,1],$ respectively $I(\subset \mathbb{R})$, and  let $g: \mathcal{D}_{1} \to \mathcal{D}_{2}$  be a conformal map such that $g(J)=I$. Let $E_{\rho}\left(\subset \mathcal{D}_{1}\right), \rho>1$, denote an ellipse with foci at $\pm 1$ and with the sum of its major and minor axes equal to $2 \rho$. If $f:\mathcal{D}_{2} \to \mathbb{C}$ is a function  such that the composition $f \circ g: \mathcal{D}_{1} \to \mathbb{C}$ is analytic inside and on $E_{\rho}$, and if $g$ is analytic inside and on $E_\sigma \supseteq E_\rho$, then the rational function \eqref{eq:rat} equipping weights \eqref{eq:weig} and interpolating $f$ between the transformed Chebyshev points $x_{i}:=g(y_{i})$, satisfies
  \begin{align*}
    | f(x)-r(x) | = \mathcal{O}\lrbrack{\rho^{-n}}
  \end{align*}
  uniformly for all $x \in[-1,1]$.
\end{theorem}

The proof of {\sc Theorem}~\ref{thm:cite1999} is applied  the  map $w : \mathcal{D}_1 \times \mathcal{D}_2 \to \mathbb{C}$ as $w(z,y)=\frac{z-y}{g(z)-g(y)}$ which is analytic for each fixed $z\in  \mathcal{D}_1$ and satisfies $\min |w(y,y)|>0$.
In this paper, however, as pointed out in Introduction, the map $w$  defined as
\begin{align*}
  w(z,y)=\frac{z-y}{g(z)-g(y)},\quad w(y,y)=\frac{2^{s/\alpha} \alpha}{Ts \lrbrack{ y+1 }^{s/\alpha-1}} \quad
  \text{for} \quad g(y)=T \lrbrack{ \frac{y+1}{2} }^{s/\alpha}
\end{align*}
has poles at $z=y=-1 \in E_\rho$ due to that $\alpha \in (0,1)$,  $s$ is a positive integer and  then $s/\alpha>1$.

In the following, we adopt the idea proposed in \cite{Guettel2012} to analyse the rational approximation \eqref{eq:rat} on a subset $[\delta, T] \subset \Omega$ for $0 < \delta_0 < \delta < T$.
\begin{lemma}\label{lem:pot1}
  Let $g: [-1,1] \to \Omega \subset \mathbb{R}$ be a monotonically increasing map, then for $z \in \mathbb{C} \backslash \Omega$, we have
  \begin{align}\label{eq:lem1_0}
    \limsup_{N \to \infty} \snorm{ \sum_{k=0}^N \frac{\lambda_k}{z-x_k} }^{1/(N+1)}
    \leq \exp\lrbrack{-U_{\hat{y}}},
  \end{align}
  where
  \begin{align}\label{eq:lem1_1}
    \exp\lrbrack{-U_{\hat{y}}} := \max_{ y_k \in [-1,1] } \exp\lrbrack{-U_{y_k}},
  \end{align}
  with
  \begin{align}\label{eq:lem1_2}
    U_{y_k} := \frac{1+y_k}{2} \int_{-1}^{y_k} \log|y_k-y| \phi(y) dy + \frac{1-y_k}{2} \int_{y_k}^{1} \log|y_k-y| \phi(y) dy,
  \end{align}
  and
  \begin{align}\label{eq:lem1_3}
    \phi(y) = \frac{1}{\pi\sqrt{1-y^2}}.
  \end{align}
  \begin{proof}
    Let $l(z) = \prod_{k=0}^N (z-x_k)$ and consider
    \begin{align*}
      \sum_{k=0}^N \snorm{ l(z) \frac{\lambda_k}{z-x_k} }
      &= \sum_{k=0}^N \prod_{j=0, j \neq k}^N \snorm{ \frac{z-x_j}{y_k-y_j} }
      = \sum_{k=0}^N \exp\lrbrack{ \sum_{j=0, j \neq k}^N \log \snorm{ \frac{z-x_j}{y_k-y_j} }}  \\
      &= \sum_{k=0}^N \exp\lrbrack{ N \sum_{j=0, j \neq k}^N \frac{\log|z-x_j|}{N} - N\sum_{j=0, j \neq k}^N \frac{\log|y_k-y_j|}{N} }.
    \end{align*}
    Then letting $N$ tend to infinity on the both sides of the above formula, we have
    \begin{align*}
      \lim_{N \to \infty} \sum_{k=0}^N \snorm{ l(z) \frac{\lambda_k}{z-x_k} }
      = \sum_{k=0}^N \exp\lrbrack{ N U(g(z)) - NU_{y_k} },
    \end{align*}
    where
    \begin{align}\label{eq:lem1_4}
      U(g(z)) := \int_{-1}^1 \log|g(z)-g(y)| \phi(y) dy < \infty,
    \end{align}
    and $U_{y_k}$ and $\phi(y)$ are given by \eqref{eq:lem1_2} and \eqref{eq:lem1_3}, respectively.
    With the same operations, it can be obtained that
    \begin{align*}
       \lim_{N \to \infty} \snorm{l(z)} = \exp((N+1)U(g(z))).
    \end{align*}
    Thus, for $z \in \mathbb{C} \backslash [-1,1]$, there holds
    \begin{align*}
      &\limsup_{N \to \infty} \snorm{ \sum_{k=0}^N \frac{\lambda_k}{z-x_k} }^{N+1}
      = \limsup_{N \to \infty} \snorm{ \left. l(z)\sum\limits_{k=0}^N \frac{\lambda_k}{z-x_k} \right/ l(z) }^{N+1}  \\
      & = \limsup_{N \to \infty} \snorm{ \left. \sum\limits_{k=0}^N \exp\lrbrack{ NU(g(z)) - NU_{y_k} } \right/ \exp((N+1)U(g(z))) }^{N+1} \\
      &= \exp\lrbrack{-U_{\hat{y}}}.
    \end{align*}
    From the definition of \eqref{eq:lem1_1}, we prove the result.
  \end{proof}
\end{lemma}
\begin{corollary}\label{col:den}
  Let $g: [\sigma,1] \to [\delta, T] \subset \mathbb{R}$ be a monotonically increasing map, then for $z \in \mathbb{C} \backslash [\delta, T]$, we have
  \begin{align*}
    \limsup_{N \to \infty} \snorm{ \sum_{k=0}^N \frac{\lambda_k}{z-x_k} }^{1/(N+1)}
    \leq \exp(-\tilde{U}_{\hat{y}}),
  \end{align*}
  where
  \begin{align*}
    \exp(-\tilde{U}_{\hat{y}}) := \max_{ y_k \in [\sigma,1] } \exp(-\tilde{U}_{y_k}),
  \end{align*}
  with
  \begin{align*}
    \tilde{U}_{y_k} := \frac{1+y_k}{2} \int_{\sigma}^{y_k} \log|y_k-y| \phi(y) dy + \frac{1-y_k}{2} \int_{y_k}^{1} \log|y_k-y| \phi(y) dy.
  \end{align*}
\end{corollary}

\begin{lemma}\label{lem:pot2}
  Let $g: [-1,1] \to \Omega \subset \mathbb{R}$ be a strictly monotonic increasing map, then for $x \in \Omega$ (or equivalently $x=g(y), y \in [-1,1]$), we have
  \begin{align}\label{eq:lem2_1}
    \liminf_{N \to \infty} \snorm{ \sum_{k=0}^N \frac{\lambda_k}{x-x_k} }^{1/(N+1)} = \exp(U_{g(y)}-U_{\tilde{y}}-U(g(y))),
  \end{align}
  where $U(g(y))$ is defined by \eqref{eq:lem1_4},
  \begin{align*}
    U_{\tilde{y}} = \max_{y_k \in [-1,1]} U_{y_k},
  \end{align*}
  and
  \begin{align}\label{eq:lem2_2}
    U_{g(y)} := \frac{1+y}{2} \int_{-1}^y \log|g(y)-g(z)| \phi(z) dz + \frac{1-y}{2} \int_{y}^{1} \log|g(y)-g(z)| \phi(z) dz.
  \end{align}
\begin{proof}
  The same computation as in {\sc Lemma}~{lem:pot1} but now for $x \in \Omega$ yields
  \begin{align*}
    \sum_{k=0}^N \snorm{ l(x)  \frac{\lambda_k}{x-x_k} }
    = \sum_{k=0}^N \exp\lrbrack{ N \sum_{j=0, j \neq k}^N \frac{\log|x-x_j|}{N} - N\sum_{j=0, j \neq k}^N \frac{\log|y_k-y_j|}{N} }.
  \end{align*}
  Taking the limit as $N \to \infty$ on both sides of the above formula yields
  \begin{align*}
    \lim_{N \to \infty} \sum_{j=0, j \neq k}^N \frac{\log|x-x_j|}{N}
    = U_{g(y)}, \mbox{ and }
    \lim_{N \to \infty} \sum_{j=0, j \neq k}^N \frac{\log|y_k-y_j|}{N} = U_{y_k},
  \end{align*}
  where $U_{g(y)}$ and $U_{y_k}$ are defined respectively in \eqref{eq:lem2_2} and \eqref{eq:lem1_2}.

  The lemma holds provided
  \begin{align*}
    \liminf_{N \to \infty} \snorm{ \sum_{k=0}^N \frac{\lambda_k}{x-x_k} }^{1/(N+1)}
    &= \liminf_{N \to \infty} \snorm{ \left. l(x)\sum\limits_{k=0}^N \frac{\lambda_k}{x-x_k} \right/ l(x) }^{1/(N+1)} \\
    &= \frac{\exp(U_{g(y)}-U_{\tilde{y}})}{\exp(U(g(y)))}.
  \end{align*}
\end{proof}
\end{lemma}
\begin{corollary}\label{col:num}
  Let $g: [\sigma,1] \to [\delta,T] \subset \mathbb{R}$ be a strictly monotonic increasing map, then for $x \in [\delta, T]$ (or equivalently $x=g(y), y \in [\sigma,1]$), we have
  \begin{align*}
    \liminf_{N \to \infty} \snorm{ \sum_{k=0}^N \frac{\lambda_k}{x-x_k} }^{1/(N+1)} = \exp(\tilde{U}_{g(y)}-\tilde{U}_{\tilde{y}}-\tilde{U}(g(y))),
  \end{align*}
  where $\tilde{U}(g(y))$ is defined as
  \begin{align*}
    \tilde{U}(g(y)) := \int_{\sigma}^1 \log|g(y)-g(w)| \phi(w) dw,
  \end{align*}
  and
  \begin{align*}
    \tilde{U}_{\tilde{y}} = \max_{y_k \in [\sigma,1]} U_{y_k},
  \end{align*}
  with
  \begin{align*}
    \tilde{U}_{g(y)} := \frac{1+y}{2} \int_{\sigma}^y \log|g(y)-g(z)| \phi(z) dz + \frac{1-y}{2} \int_{y}^{1} \log|g(y)-g(z)| \phi(z) dz.
  \end{align*}
\end{corollary}

We  begin investigating the asymptotic convergence of the linear barycentric rational interpolation for function $f(x)$ having a singularity at the origin. Pointed in \cite{Guettel2012}, it is sufficient to investigate interpolants, of ``prototype functions'' $h(x,z) = 1/(z-x)$ with a simple poles $z \in \mathbb{C} \backslash (0, T]$. An explicit formula for the polynomial interpolants $p$ of such a particular function is
\begin{align*}
  p_i(x) = \left. \lrbrack{1-\dfrac{x-x_i}{z-x_i}} \right/ (z-x).
\end{align*}
It can be verified that this is indeed a constant satisfying $p_i(x_i) = h(x_i,z)$. Hence the rational interpolant \eqref{eq:rat} of $h$, which we denote by $r[h](x)$, is
\begin{align*}
  r[h](x) = \frac{1}{z-x} \dfrac{\sum\limits_{i=0}^N \dfrac{\lambda_i}{x-x_i} \lrbrack{1-\dfrac{x-x_i}{z-x_i}} }{\sum\limits_{i=0}^N \dfrac{\lambda_i}{x-x_i} }
  = \frac{1}{z-x} \lrbrack{ 1- \left. \sum\limits_{i=0}^N \dfrac{\lambda_i}{z-x_i} \right/ \sum\limits_{i=0}^N \dfrac{\lambda_i}{x-x_i} },
\end{align*}
and, therefore,
\begin{align*}
  h(x,z) - r[h](x) = \frac{1}{z-x} \left. \sum\limits_{i=0}^N \dfrac{\lambda_i}{z-x_i} \right/ \sum\limits_{i=0}^N \dfrac{\lambda_i}{x-x_i}.
\end{align*}
Let $V(z)$ be the ``potential function'' defined as
\begin{align*}
  V(g(z)) :=
  \begin{dcases}
    -\tilde{U}_{\hat{y}}, & z \in \mathbb{C} \backslash [-1,1],  \\
    \tilde{U}_{g(z)}-\tilde{U}_{\tilde{y}}-\tilde{U}(g(z)), & z \in [-1,1].
  \end{dcases}
\end{align*}
Combining {\sc Corollary} \ref{col:den} and {\sc Corollary} \ref{col:num}, and considering the monotonicity of the exponential function, we can deduce that
\begin{align*}
  \limsup_{N \to \infty} \snorm{ h(x,z)-r[h](x) }^{1/(N+1)} \leq \exp\lrbrack{V(g(z)) - V(g(y))}.
\end{align*}

Next, let us turn to analyse the uniform convergence of the linear barycentric rational interpolation \eqref{eq:rat} over the whole interval $\Omega$.
For this purpose, let us define the contours
\begin{align}\label{eq:cr}
  \mathcal{C}_R := \lrBbrack{z \in \mathbb{C}: R:= \dfrac{\exp\lrbrack{V(g(z))}}{\min_{y \in [\sigma,1]} \exp\lrbrack{V(g(y))}} },
\end{align}
which can be seen as levels of convergence with rate  at least $R$ for every points $x \in \Omega$.

For function $f$ which is analytic on the region in complex plane containing $\Omega \backslash [0, \delta)$ for $0 < \delta_0 < \delta < T$, it can be represented by Cauchy integral formula
\begin{align*}
  f(x) = \frac{1}{2\pi i} \int_\mathcal{C} \frac{f(z)}{z-x} dz = \frac{1}{2\pi i} \int_\mathcal{C} f(z) h(x,z) dz,
\end{align*}
which, together with the rational approximation $r[h]$, implies that
\begin{align*}
  r(x) = \frac{1}{2\pi i} \int_\mathcal{C} \frac{f(z)}{z-x} \lrbrack{ 1- \left. \sum\limits_{i=0}^N \dfrac{\lambda_i}{z-x_i} \right/ \sum\limits_{i=0}^N \dfrac{\lambda_i}{x-x_i} } dz
\end{align*}
is the rational interpolant for $f$. Therefore, the interpolation error is
\begin{align*}
  f(x)-r(x) = \frac{1}{2\pi i} \int_\mathcal{C} \frac{f(z)}{z-x} \left. \sum\limits_{i=0}^N \dfrac{\lambda_i}{z-x_i} \right/ \sum\limits_{i=0}^N \dfrac{\lambda_i}{x-x_i} dz.
\end{align*}
Then, for $\delta>\delta_0>0$, we have
\begin{align*}
  \norm{f-r}_\infty \leq D \left. \sup_{z \in \mathcal{C}} \snorm{ \sum\limits_{i=0}^N \dfrac{\lambda_i}{z-x_i} } \right/ \min_{x\in [\delta,T]} \sum\limits_{i=0}^N \dfrac{\lambda_i}{x-x_i},
\end{align*}
where $D = \frac{\mbox{length}(\mathcal{C}) \max_{z \in \mathcal{C}}|f(z)|}{2\pi \mbox{dist}([\delta,T], \mathcal{C})}$ is a constant independent of $N$. Thus, we can obtain the convergence rate of \eqref{eq:rat} for $f$ which is expressed in the following theorem.
\begin{theorem}\label{thm:conv}
  Let $f$ be a function analytic in an open region containing $[\delta, T] \subset \Omega$ for $0 < \delta_0 < \delta$, and let $R$ be the smallest number such that $f$ is analytic in the interior of $\mathcal{C}_R$ defined in \eqref{eq:cr}. Then the rational interpolants $r$ defined by \eqref{eq:rat} with map \eqref{eq:refpts} or \eqref{eq:refpts2} satisfy
  \begin{align*}
    \limsup_{N \to \infty} \norm{ f-r }_\infty^{1/N} \leq R.
  \end{align*}
\end{theorem}

\begin{figure}[th]
  \centering
  \includegraphics[width=0.9\textwidth]{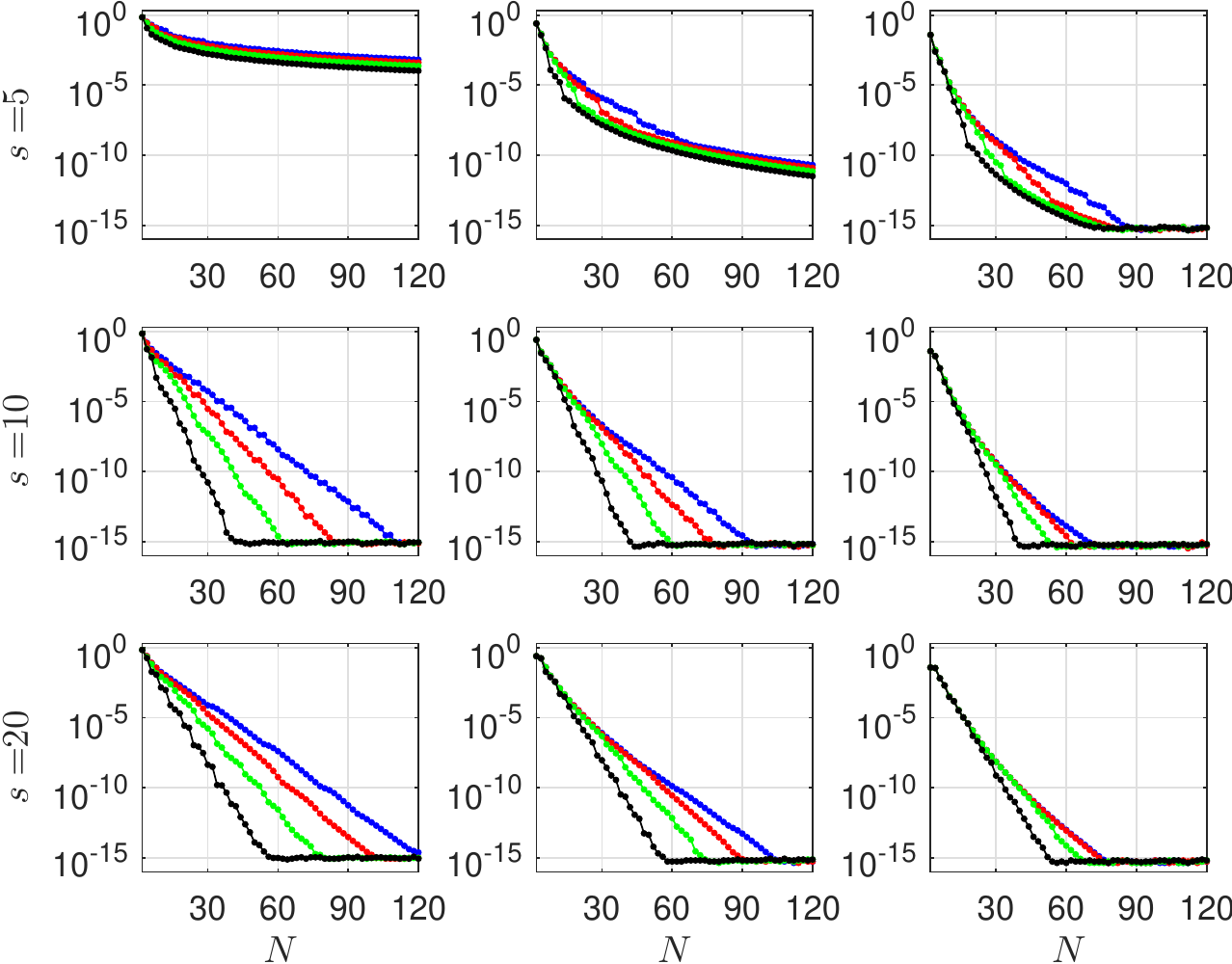}
  \caption{Errors of \eqref{eq:rat} on subset $[\delta, 1] \subset [0,1]$ against $N$ varying evenly from $2$ to $120$ in maximum error in \texttt{{\rm x0 = logspace(log10(del), 0, 10000)}} for function $x^\alpha$ with $\alpha=0.1$ (first column),  $\alpha=0.5$ (second column), $\alpha=0.9$ (third column), and $s=5$ (first row), $s=10$ (second row), $s=20$ (third row) and $\delta = 10^{-7}$ (blue), $10^{-5}$ (red), $10^{-3}$ (green), $10^{-1}$ (black).}
  \label{fig:ex1_abs_del}
\end{figure}

{\sc Theorem}~\ref{thm:conv} implies the convergence rate of \eqref{eq:rat} is related to $R$ introduced in \eqref{eq:cr}.
Undoubtedly, the exponentially convergent rate also depends on $\delta$ and $s$, shown in {\sc Fig.}~\ref{fig:ex1_abs_del}, in which function $x^\alpha$ is considered on $[\delta, 1] \subset [0,1]$ with different values of $\alpha$ and $s$. The maximum error is calculated in set \texttt{x0 = logspace(log10(del), 0, 10000)}, where \texttt{del} denotes $\delta$.

\begin{figure}[th]
  \centering
  \includegraphics[width=0.45\textwidth]{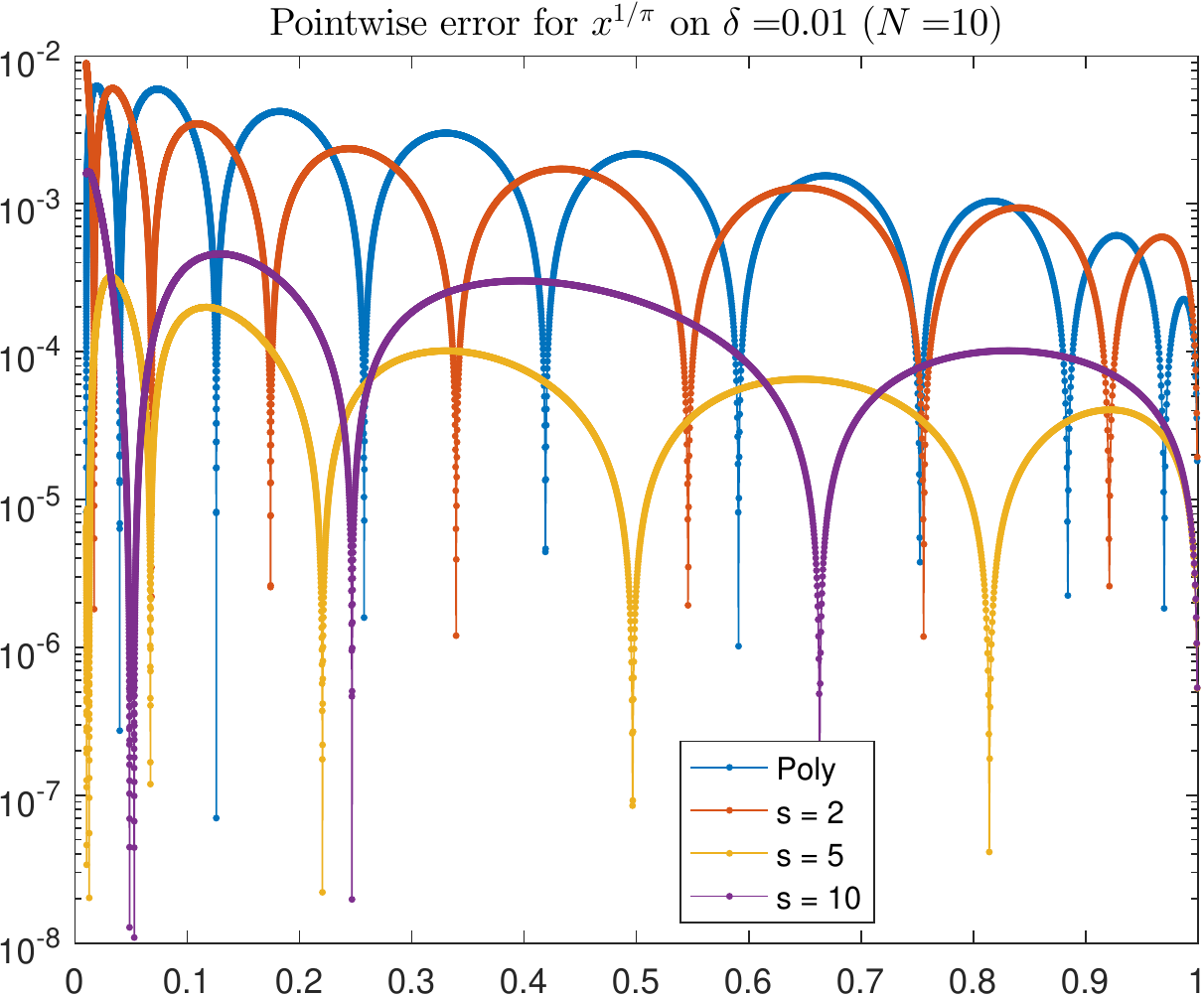}
  \includegraphics[width=0.45\textwidth]{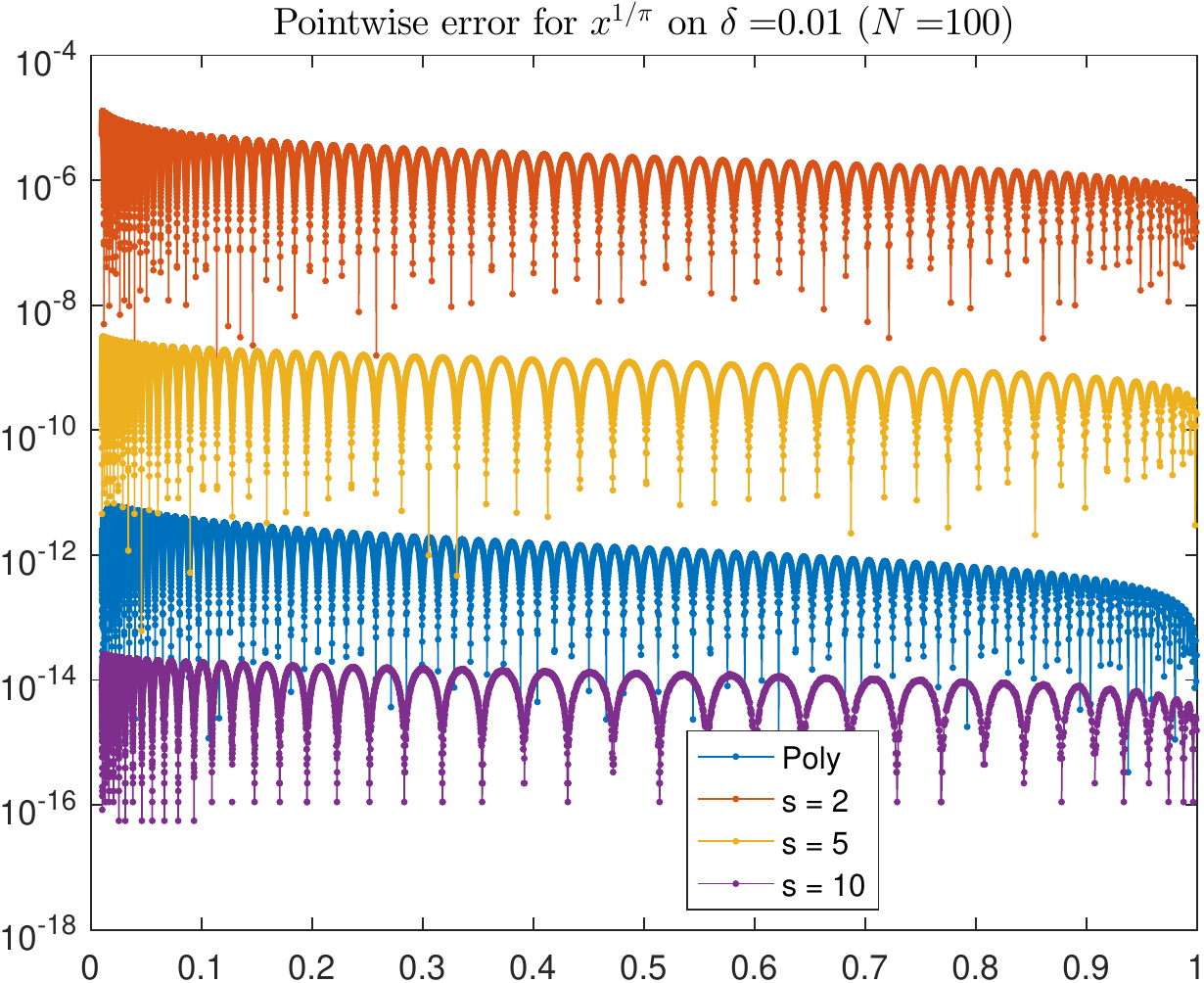}  \\
  \includegraphics[width=0.45\textwidth]{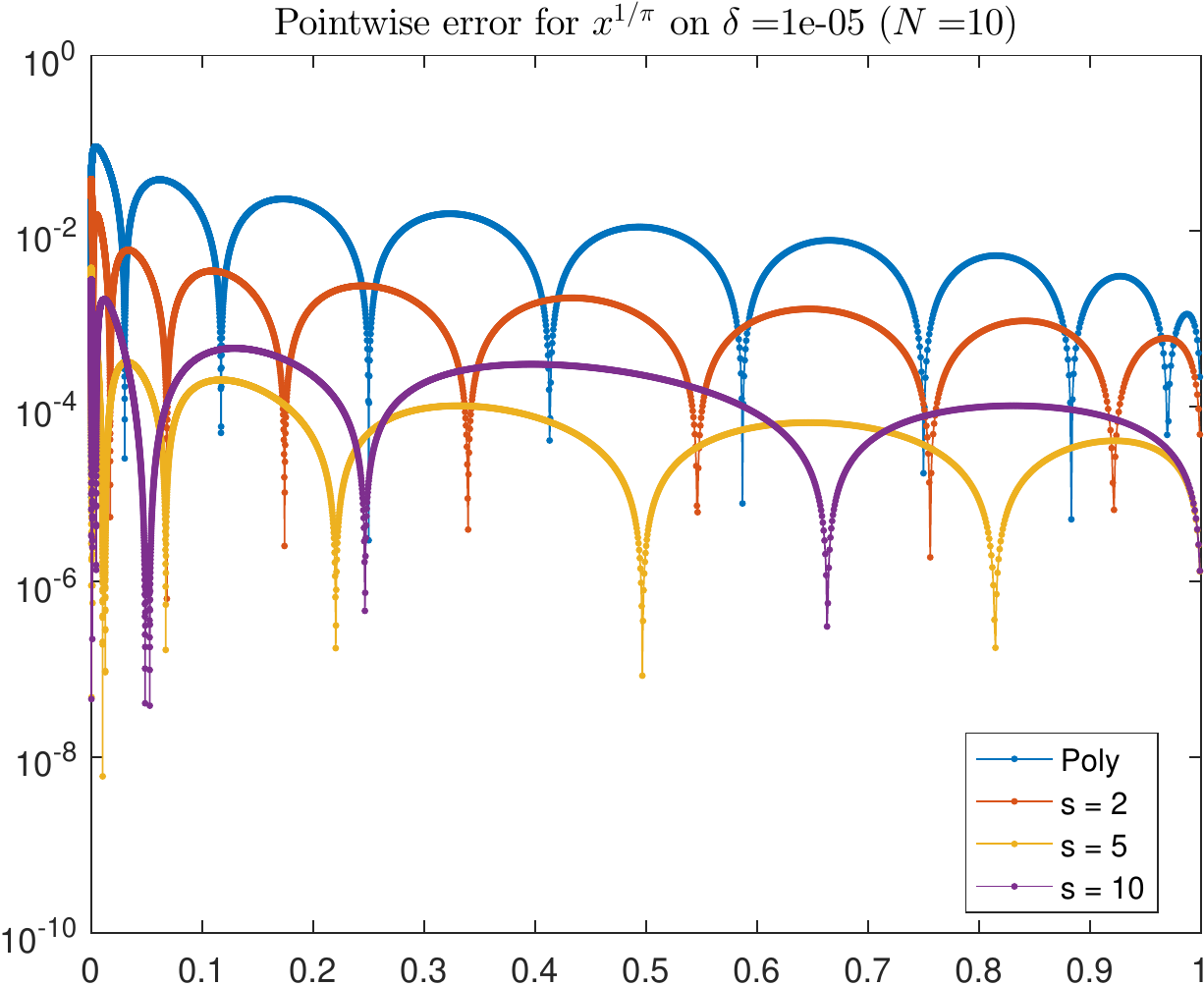}
  \includegraphics[width=0.45\textwidth]{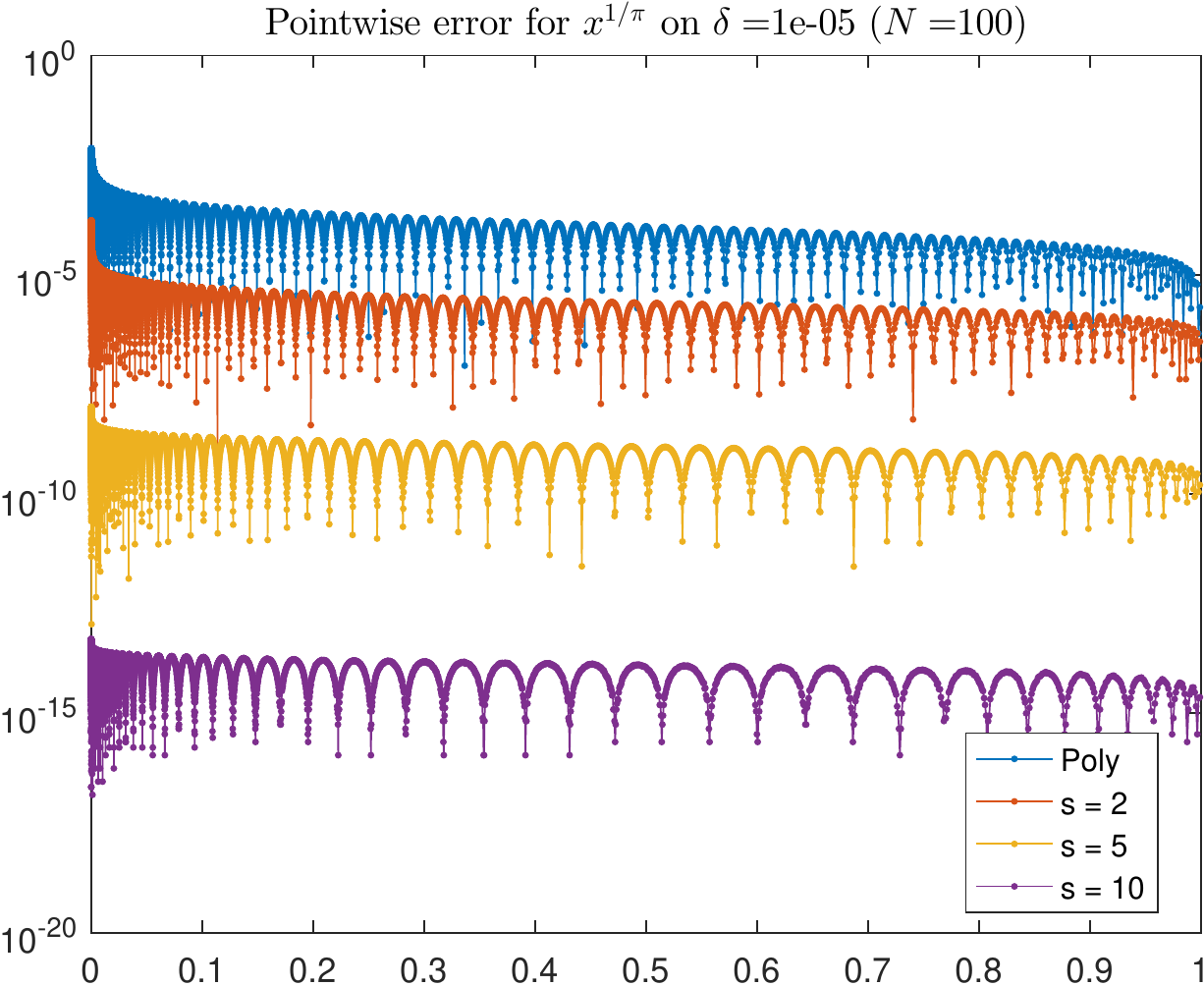}
  \caption{Absolute errors of the polynomial interpolation  at \eqref{SCheb}  and the rational interpolation \eqref{eq:rat} associated with Legendre weights and scaled points  on subset $[\delta, 1] \subset [0,1]$  for function $x^\alpha$ with $\alpha=1/\pi$,  $\delta=10^{-2}$ (first row) and  $\delta=10^{-5}$ (second row), respectively. In each plots, different values of $s$ in \eqref{eq:refpts2} are considered with $N = 10$ (first column) and $N = 100$ (second column).}
  \label{fig:polSca1}
\end{figure}

It is of particular interest to notice that for approximation of  functions near the singular point  $x=0$, for example, $f(x)=x^{\alpha}$ in $[\delta,1]$ with $0<\delta\ll 1$,
the the proposed rational approximation \eqref{eq:rat}  performs much better than the corresponding polynomial interpolation at the shifted Chebyshev points of second kind
\begin{equation}\label{SCheb}
x_i=\frac{1+\delta}{2}+\frac{1-\delta}{2}\cos(i\pi/N),\quad i=0,1,\ldots,N.
\end{equation}
In {\sc Fig.}~\ref{fig:polSca1}, we simulate the rational approximation \eqref{eq:rat} associated with Legendre weights and scaled points by \eqref{eq:refpts2}. The simulated results with these weights and points are similar to \eqref{eq:rat} associated with the Chebyshev case, which illustrates again the proposed method is also efficient in Jacobi-Gauss-Lobatto situations.

\section{Conclusions}\label{sec:condis}

This paper is concerned with the linear barycentric rational interpolant  coupled with strictly  monotonic increasing maps \eqref{eq:refpts} or \eqref{eq:refpts2} for $x^\alpha$-type functions ($\alpha \in (0,1)$), or \eqref{eq:refpts_log} for logarithmic type functions, which is easily implemented and just take $O(N)$ operations.
Various numerical experiments also illustrate it's  well-condition and accuracy.


The function approximation associated with the roots or extrema of Jacobi orthogonal polynomial is widely popular in approximation theory. 
Numerical examples show that the linear barycentric rational interpolant  with the scaled map \eqref{eq:refpts} for $x^\alpha$ while \eqref{eq:refpts_log} for the  logarithmic  function adopted the Jacobi-Gauss, Jacobi-Gauss-Lobatto or Jacobi-Gauss-Radau points with respect to the corrsponding simplified weights  \cite{Wang2014} is efficient too.

However, there are still many problems should be considered further such as  the convergence rates and distribution of the poles.

Based on   a series of  numerical experiments, we give the following hypothesis on the convergence of \eqref{eq:rat} for $|x|^\alpha$.
\begin{hypothesis}\label{conj1}
  Assume $r(x)$ is the linear barycentric rational interpolation given by \eqref{eq:rat} with \eqref{eq:refpts} to $|x|^\alpha$, then for $\alpha \in (0,1]$, the interpolation error satisfies
  \begin{align*}
    \norm{ |x|^\alpha - r(x)}_\infty =
    \begin{dcases}
      O(N^{-2s}), & s < s_0, \\
      O(\exp(-C_{\alpha,s} \sqrt{N})), & s > s_1,
    \end{dcases}
  \end{align*}
  for some positive numbers $0 < s_0 < s_1$ and $C_{\alpha,s}$ depending only on $\alpha$ and $s$.
  Equivalent result holds for rational interpolation \eqref{eq:rat} with \eqref{eq:refpts2}:
  \begin{align*}
    \norm{ |x|^\alpha - r(x)}_\infty =
    \begin{dcases}
      O(N^{-2\alpha s}), & s < s_0, \\
      O(\exp(-C_{\alpha,s} \sqrt{N})), & s > s_1.
    \end{dcases}
  \end{align*}
\end{hypothesis}

\bibliographystyle{siamplain}

\end{document}